\documentclass[reqno,12pt,a4paper]{amsart}

\usepackage{mathrsfs}
\usepackage{amssymb}
\usepackage{amsfonts}
\usepackage{latexsym}
\usepackage{amsthm}
\usepackage{graphicx}

\def\lb{\label}

\begin{document}


\renewcommand{\theequation}{\arabic{section}.\arabic{equation}}
\theoremstyle{plain}
\newtheorem{theorem}{\bf Theorem}[section]
\newtheorem{lemma}[theorem]{\bf Lemma}
\newtheorem{corollary}[theorem]{\bf Corollary}
\newtheorem{proposition}[theorem]{\bf Proposition}
\newtheorem{definition}[theorem]{\bf Definition}

\newtheorem{remark}[theorem]{\bf Remark}

\def\a{\alpha}  \def\cA{{\mathcal A}}     \def\bA{{\bf A}}  \def\mA{{\mathscr A}}
\def\b{\beta}   \def\cB{{\mathcal B}}     \def\bB{{\bf B}}  \def\mB{{\mathscr B}}
\def\g{\gamma}  \def\cC{{\mathcal C}}     \def\bC{{\bf C}}  \def\mC{{\mathscr C}}
\def\G{\Gamma}  \def\cD{{\mathcal D}}     \def\bD{{\bf D}}  \def\mD{{\mathscr D}}
\def\d{\delta}  \def\cE{{\mathcal E}}     \def\bE{{\bf E}}  \def\mE{{\mathscr E}}
\def\D{\Delta}  \def\cF{{\mathcal F}}     \def\bF{{\bf F}}  \def\mF{{\mathscr F}}
\def\c{\chi}    \def\cG{{\mathcal G}}     \def\bG{{\bf G}}  \def\mG{{\mathscr G}}
\def\z{\zeta}   \def\cH{{\mathcal H}}     \def\bH{{\bf H}}  \def\mH{{\mathscr H}}
\def\e{\eta}    \def\cI{{\mathcal I}}     \def\bI{{\bf I}}  \def\mI{{\mathscr I}}
\def\p{\psi}    \def\cJ{{\mathcal J}}     \def\bJ{{\bf J}}  \def\mJ{{\mathscr J}}
\def\vT{\Theta} \def\cK{{\mathcal K}}     \def\bK{{\bf K}}  \def\mK{{\mathscr K}}
\def\k{\kappa}  \def\cL{{\mathcal L}}     \def\bL{{\bf L}}  \def\mL{{\mathscr L}}
\def\l{\lambda} \def\cM{{\mathcal M}}     \def\bM{{\bf M}}  \def\mM{{\mathscr M}}
\def\L{\Lambda} \def\cN{{\mathcal N}}     \def\bN{{\bf N}}  \def\mN{{\mathscr N}}
\def\m{\mu}     \def\cO{{\mathcal O}}     \def\bO{{\bf O}}  \def\mO{{\mathscr O}}
\def\n{\nu}     \def\cP{{\mathcal P}}     \def\bP{{\bf P}}  \def\mP{{\mathscr P}}
\def\r{\rho}    \def\cQ{{\mathcal Q}}     \def\bQ{{\bf Q}}  \def\mQ{{\mathscr Q}}
\def\s{\sigma}  \def\cR{{\mathcal R}}     \def\bR{{\bf R}}  \def\mR{{\mathscr R}}
\def\S{\Sigma}  \def\cS{{\mathcal S}}     \def\bS{{\bf S}}  \def\mS{{\mathscr S}}
\def\t{\tau}    \def\cT{{\mathcal T}}     \def\bT{{\bf T}}  \def\mT{{\mathscr T}}
\def\f{\phi}    \def\cU{{\mathcal U}}     \def\bU{{\bf U}}  \def\mU{{\mathscr U}}
\def\F{\Phi}    \def\cV{{\mathcal V}}     \def\bV{{\bf V}}  \def\mV{{\mathscr V}}
\def\P{\Psi}    \def\cW{{\mathcal W}}     \def\bW{{\bf W}}  \def\mW{{\mathscr W}}
\def\o{\omega}  \def\cX{{\mathcal X}}     \def\bX{{\bf X}}  \def\mX{{\mathscr X}}
\def\x{\xi}     \def\cY{{\mathcal Y}}     \def\bY{{\bf Y}}  \def\mY{{\mathscr Y}}
\def\X{\Xi}     \def\cZ{{\mathcal Z}}     \def\bZ{{\bf Z}}  \def\mZ{{\mathscr Z}}
\def\be{{\bf e}}
\def\bv{{\bf v}} \def\bu{{\bf u}}
\def\Om{\Omega}
\def\bbD{\pmb \Delta}
\def\mm{\mathrm m}
\def\mn{\mathrm n}

\newcommand{\mc}{\mathscr {c}}

\newcommand{\gA}{\mathfrak{A}}          \newcommand{\ga}{\mathfrak{a}}
\newcommand{\gB}{\mathfrak{B}}          \newcommand{\gb}{\mathfrak{b}}
\newcommand{\gC}{\mathfrak{C}}          \newcommand{\gc}{\mathfrak{c}}
\newcommand{\gD}{\mathfrak{D}}          \newcommand{\gd}{\mathfrak{d}}
\newcommand{\gE}{\mathfrak{E}}
\newcommand{\gF}{\mathfrak{F}}           \newcommand{\gf}{\mathfrak{f}}
\newcommand{\gG}{\mathfrak{G}}           
\newcommand{\gH}{\mathfrak{H}}           \newcommand{\gh}{\mathfrak{h}}
\newcommand{\gI}{\mathfrak{I}}           \newcommand{\gi}{\mathfrak{i}}
\newcommand{\gJ}{\mathfrak{J}}           \newcommand{\gj}{\mathfrak{j}}
\newcommand{\gK}{\mathfrak{K}}            \newcommand{\gk}{\mathfrak{k}}
\newcommand{\gL}{\mathfrak{L}}            \newcommand{\gl}{\mathfrak{l}}
\newcommand{\gM}{\mathfrak{M}}            \newcommand{\gm}{\mathfrak{m}}
\newcommand{\gN}{\mathfrak{N}}            \newcommand{\gn}{\mathfrak{n}}
\newcommand{\gO}{\mathfrak{O}}
\newcommand{\gP}{\mathfrak{P}}             \newcommand{\gp}{\mathfrak{p}}
\newcommand{\gQ}{\mathfrak{Q}}             \newcommand{\gq}{\mathfrak{q}}
\newcommand{\gR}{\mathfrak{R}}             \newcommand{\gr}{\mathfrak{r}}
\newcommand{\gS}{\mathfrak{S}}              \newcommand{\gs}{\mathfrak{s}}
\newcommand{\gT}{\mathfrak{T}}             \newcommand{\gt}{\mathfrak{t}}
\newcommand{\gU}{\mathfrak{U}}             \newcommand{\gu}{\mathfrak{u}}
\newcommand{\gV}{\mathfrak{V}}             \newcommand{\gv}{\mathfrak{v}}
\newcommand{\gW}{\mathfrak{W}}             \newcommand{\gw}{\mathfrak{w}}
\newcommand{\gX}{\mathfrak{X}}               \newcommand{\gx}{\mathfrak{x}}
\newcommand{\gY}{\mathfrak{Y}}              \newcommand{\gy}{\mathfrak{y}}
\newcommand{\gZ}{\mathfrak{Z}}             \newcommand{\gz}{\mathfrak{z}}

\def\ve{\varepsilon}   \def\vt{\vartheta}    \def\vp{\varphi}    \def\vk{\varkappa}

\def\A{{\mathbb A}} \def\B{{\mathbb B}} \def\C{{\mathbb C}}
\def\dD{{\mathbb D}} \def\E{{\mathbb E}} \def\dF{{\mathbb F}} \def\dG{{\mathbb G}} \def\H{{\mathbb H}}\def\I{{\mathbb I}} \def\J{{\mathbb J}} \def\K{{\mathbb K}} \def\dL{{\mathbb L}}\def\M{{\mathbb M}} \def\N{{\mathbb N}} \def\O{{\mathbb O}} \def\dP{{\mathbb P}} \def\R{{\mathbb R}}\def\S{{\mathbb S}} \def\T{{\mathbb T}} \def\U{{\mathbb U}} \def\V{{\mathbb V}}\def\W{{\mathbb W}} \def\X{{\mathbb X}} \def\Y{{\mathbb Y}} \def\Z{{\mathbb Z}}


\def\la{\leftarrow}              \def\ra{\rightarrow}            \def\Ra{\Rightarrow}
\def\ua{\uparrow}                \def\da{\downarrow}
\def\lra{\leftrightarrow}        \def\Lra{\Leftrightarrow}


\def\lt{\biggl}                  \def\rt{\biggr}
\def\ol{\overline}               \def\wt{\widetilde}
\def\ul{\underline}
\def\no{\noindent}


\let\ge\geqslant                 \let\le\leqslant
\def\lan{\langle}                \def\ran{\rangle}
\def\/{\over}                    \def\iy{\infty}
\def\sm{\setminus}               \def\es{\emptyset}
\def\ss{\subset}                 \def\ts{\times}
\def\pa{\partial}                \def\os{\oplus}
\def\om{\ominus}                 \def\ev{\equiv}
\def\iint{\int\!\!\!\int}        \def\iintt{\mathop{\int\!\!\int\!\!\dots\!\!\int}\limits}
\def\el2{\ell^{\,2}}             \def\1{1\!\!1}
\def\sh{\sharp}
\def\wh{\widehat}
\def\bs{\backslash}
\def\intl{\int\limits}

\def\na{\mathop{\mathrm{\nabla}}\nolimits}
\def\sh{\mathop{\mathrm{sh}}\nolimits}
\def\ch{\mathop{\mathrm{ch}}\nolimits}
\def\where{\mathop{\mathrm{where}}\nolimits}
\def\all{\mathop{\mathrm{all}}\nolimits}
\def\as{\mathop{\mathrm{as}}\nolimits}
\def\Area{\mathop{\mathrm{Area}}\nolimits}
\def\arg{\mathop{\mathrm{arg}}\nolimits}
\def\const{\mathop{\mathrm{const}}\nolimits}
\def\det{\mathop{\mathrm{det}}\nolimits}
\def\diag{\mathop{\mathrm{diag}}\nolimits}
\def\diam{\mathop{\mathrm{diam}}\nolimits}
\def\dim{\mathop{\mathrm{dim}}\nolimits}
\def\dist{\mathop{\mathrm{dist}}\nolimits}
\def\Im{\mathop{\mathrm{Im}}\nolimits}
\def\Iso{\mathop{\mathrm{Iso}}\nolimits}
\def\Ker{\mathop{\mathrm{Ker}}\nolimits}
\def\Lip{\mathop{\mathrm{Lip}}\nolimits}
\def\rank{\mathop{\mathrm{rank}}\limits}
\def\Ran{\mathop{\mathrm{Ran}}\nolimits}
\def\Re{\mathop{\mathrm{Re}}\nolimits}
\def\Res{\mathop{\mathrm{Res}}\nolimits}
\def\res{\mathop{\mathrm{res}}\limits}
\def\sign{\mathop{\mathrm{sign}}\nolimits}
\def\span{\mathop{\mathrm{span}}\nolimits}
\def\supp{\mathop{\mathrm{supp}}\nolimits}
\def\Tr{\mathop{\mathrm{Tr}}\nolimits}
\def\BBox{\hspace{1mm}\vrule height6pt width5.5pt depth0pt \hspace{6pt}}


\newcommand\nh[2]{\widehat{#1}\vphantom{#1}^{(#2)}}
\def\dia{\diamond}

\def\Oplus{\bigoplus\nolimits}



\def\qqq{\qquad}
\def\qq{\quad}
\let\ge\geqslant
\let\le\leqslant
\let\geq\geqslant
\let\leq\leqslant
\newcommand{\ca}{\begin{cases}}
\newcommand{\ac}{\end{cases}}
\newcommand{\ma}{\begin{pmatrix}}
\newcommand{\am}{\end{pmatrix}}
\renewcommand{\[}{\begin{equation}}
\renewcommand{\]}{\end{equation}}
\def\eq{\begin{equation}}
\def\qe{\end{equation}}
\def\[{\begin{equation}}
\def\bu{\bullet}

\title[Spectral estimates for Schr\"odinger operators]
{Spectral estimates for Schr\"odinger operators on periodic discrete graphs}

\date{\today}
\author[Evgeny Korotyaev]{Evgeny Korotyaev}
\address{Department of Mathematical Analysis, Saint-Petersburg State University, Universitetskaya nab. 7/9, St. Petersburg, 199034, Russia,
\ korotyaev@gmail.com, \
e.korotyaev@spbu.ru,}
\author[Natalia Saburova]{Natalia Saburova}
\address{Department of Mathematical Analysis, Algebra and Geometry, Northern (Arctic) Federal University, Severnaya Dvina emb. 17, Arkhangelsk, 163002, Russia,
 \ n.saburova@gmail.com, \ n.saburova@narfu.ru}


\begin{abstract}
We consider normalized Laplacians and their perturbations by periodic potentials
(Schr\"odinger operators) on periodic discrete graphs. The spectrum of the operators
consists of an absolutely continuous part, which is a union of a finite number of non-degenerate
bands, and a finite number of flat bands, i.e., eigenvalues of infinite multiplicity.
We obtain estimates of the Lebesgue measure of the spectrum in terms of geometric
parameters of the graphs and show that they become identities for some class of graphs.
We determine two-sided estimates on the lengths of the first spectral bands and on the
effective masses at the bottom of the spectrum of the Laplace and Schr\"odinger operators.
In particular, these estimates yield that the first spectral band of the Schr\"odinger
operators is non-degenerate.
\end{abstract}

\maketitle

\section {\lb{Sec1}Introduction}
\setcounter{equation}{0}

Laplace operators on graphs have a lot of applications in physics and chemistry,
see, e.g., Section 7.6 in \cite{BK13}, Chapter 8 in \cite{CDS95}, the survey
\cite{Ku02} and references therein. There are two types of graphs: discrete
and metric. Laplace operators are defined on each of them.

$\bu$ A discrete graph: the difference Laplacian acts on the space of functions defined
on the vertex set of the graph. Here, vertices play the main role, and edges are
considered as relations between graph vertices. There are some definitions of the discrete
Laplace operators: normalized, combinatorial, weighted Laplacians, see, e.g.,
\cite{BK13,Ch97,HN09,MW89,S13}. The spectrum of the discrete Laplace operator on
periodic graphs consists of a finite number of bands separated by gaps.

$\bu$ A metric graph: the graph is considered as a continuous (metric)
space consisting of edges (one-dimensional spaces) connecting graph vertices.
Here the Laplacian acts on functions defined along each edge of the graph and
satisfying special boundary conditions at the vertices, which guarantee the
self-adjointness of the operator in the corresponding $L^2$-space, see, e.g.,
\cite{BK13,KoS99,P12}. The spectrum of the metric Laplace operator on
periodic graphs covers the positive semiaxis without an infinite number of gaps.

In the paper \cite{C97} an explicit relation between the spectra of the normalized
Laplacian on a discrete graph and its counterpart on the corresponding metric graph
with all edges having equal lengths is obtained. Later in the paper \cite{KS15} the eigenfunctions
of the continuous spectrum of the Laplacian on a periodic metric graph are expressed
in terms of the eigenfunctions of the continuous spectrum of the normalized Laplacian
on the corresponding discrete graph. Thus, an explicit relation between the Laplacian on
a metric graph and the normalized Laplacian on the corresponding discrete graph is
determined. Thereby, the study of the Laplacian spectrum on a metric graph with all edges of
equal lengths is reduced to the study of the spectrum of the normalized Laplacian on a discrete graph.

In this paper we consider the normalized Laplacian and its perturbations by periodic
potentials (Schr\"odinger operators) on periodic discrete graphs. We do not assume
the graph to be embedded into a Euclidean space. But in many applications
such a natural embedding exists. For instance, the tight-binding approximation is commonly
used to describe the electronic properties of real crystals
(see, e.g., \cite{A76}). A crystalline structure is modeled
by a discrete graph consisting of vertices representing
positions of atoms and edges representing chemical bonds of
atoms, by ignoring the physical characters of atoms, which
may differ from one another. The model gives
good qualitative results in many cases. Under this approach
a simple geometric model is a graph embedded into $\R^d$ ($d=2,3$) in such a way
that it is invariant with respect to the shifts by integer vectors
$m\in\Z^d$.

It is known that the spectrum of the Schr\"odinger operators with
periodic potentials on periodic discrete graphs consists
of an absolutely continuous part and a finite number of flat bands
(i.e., eigenvalues of infinite multiplicity). The absolutely continuous spectrum is a union of a finite number
of non-degenerate spectral bands separated by gaps. Here we have a
well-known problem: to estimate lengths of the spectral bands and gaps in terms
of graph parameters and potentials. In the case of the Schr\"odinger
operator $-\D+Q$ with a periodic potential $Q$ in $\R^d$ the spectrum is absolutely
continuous \cite{T73} and consists of a finite number of non-degenerate
spectral bands separated by gaps \cite{Sk87}. There are no flat bands.
We note that in the case of $d=1$ there are two-sided estimates of $L^2$-norms of
potentials in terms of
gap lengths, see \cite{K98}, \cite{K03}, and estimates of the variation
of lengths of all bands in terms of $L^2$-norms of potentials \cite{K00}.
We do not know other estimates.

\smallskip

We describe the main goals of the paper:

\noindent{\rm 1) \it to obtain two-sided estimates on the lengths of the first spectral
bands and on the effective masses at the bottom of the spectrum of the normalized Laplace
and Schr\"odinger operators on periodic discrete graphs, and to prove that
the first spectral band of the Schr\"odinger operators is non-degenerate{\rm;}

\noindent{\rm 2)} to estimate the Lebesgue measure of the spectrum and the sum
of gap lengths in terms of geometric parameters of the graphs and potentials{\rm;}

\noindent{\rm 3)} to show that the obtained estimates of the Lebesgue measure of the
spectrum are sharp, i.e., there exist periodic graphs for which these estimates become
identities{\rm;}

\noindent{\rm 4)} to describe the possible number and positions of flat bands
for normalized Laplacians on some classes of periodic graphs{\rm;} to construct a graph,
for which the number of flat bands of the Laplacian is maximal.}

The results of this paper were used essentially  in \cite{KS16} for obtaining
spectral estimates for Laplacians on metric graphs.

\subsection{Schr\"odinger operators on periodic graphs} Let $\G=(V,\cE)$
be a connected infinite graph, possibly having loops and multiple edges. Here $V$
is the set of its vertices and $\cE$ is the set of its unoriented edges.
It is convenient to assume that each unoriented edge of the graph corresponds to
two oppositely directed edges. We denote the set of all oriented edges
of the graph $\G$ by $\cA$. An edge starting at a vertex $u\in V$ and ending
at a vertex $v\in V$ will be denoted as the ordered pair $(u,v)\in\cA$ and
is said to be \emph{incident} to the vertices $u$ and $v$. Vertices $u,v\in V$
will be called \emph{adjacent} and denoted by $u\sim v$, if $(u,v)\in\cA$.
The inverse edge of $\be=(u,v)\in\cA$ will be denoted by $\ul\be=(v,u)$.
The \emph{degree} ${\vk}_v$ of the vertex $v\in V$ is the number of all edges in $\cA$ starting at $v$.

Throughout, we consider a \emph{locally finite $\Z^d$-periodic graph} $\G$, $d\geq2$,
i.e., a graph satisfying the following conditions:

1) {\it the graph $\G$ is equipped with a free action of the abelian group $\Z^d;$

{\rm2)} the quotient graph $\G_*=\G/{\Z}^d$ is finite.}

The quotient graph $\G_*=\G/{\Z}^d$ is also called the \emph{fundamental graph}
of the periodic graph $\G$. If $\G$ is embedded into the space $\R^d$,
then the fundamental graph $\G_*$ is a graph on the $d$-dimensional torus $\R^d/\Z^d$.
The fundamental graph $\G_*=(V_*,\cE_*)$ has the vertex set $V_*=V/\Z^d$,
the set $\cE_*=\cE/\Z^d$ of unoriented edges and the set $\cA_*=\cA/\Z^d$
of oriented edges.

\medskip

We consider the Hilbert space $\ell^2(V)$ of all square summable functions $f:V\to \C$
equipped with the norm
$$
\|f\|^2_{\ell^2(V)}=\sum_{v\in V}|f(v)|^2<\infty.
$$
The normalized Laplacian (the Laplace operator) $\D$ acting on $\ell^2(V)$ is defined by
\begin{equation}
\lb{DOL}
\big(\D f\big)(v)=f(v)-\sum\limits_{(v,u)\in\cA}
\frac1{\sqrt{\vk_v\vk_u}}\,f(u), \quad f\in \ell^2(V), \quad v\in V,
\end{equation}
where $\vk_v$ is the degree of the vertex $v\in V$. The sum in \eqref{DOL}
is taken over all edges in $\cA$ starting at the vertex $v$. It is known
(see \cite{MW89}) that the normalized Laplacian $\D$ is a bounded self-adjoint
operator on $\ell^2(V)$ and its spectrum $\s(\D)$ is a closed
subset of the segment $[0,2]$, containing the point 0, i.e.,
\begin{equation}
\lb{mp}
0\in\s(\D)\ss[0,2].
\end{equation}

We consider the Schr\"odinger operator $H=\D+Q$ acting on $\ell^2(V)$.
Suppose that the potential $Q$ is real valued and $\Z^d$-periodic, i.e., it satisfies
$$
Q(v+m)=Q(v), \quad  \forall\, (v,m)\in V\ts\Z^d,
$$
where $v+m$ denotes the action of $m\in\Z^d$ on $v\in V$.

\medskip

\no\textbf{Remark.} There are other definitions of Laplacians on graphs, see \cite{MW89}.
For example, the combinatorial Laplacian $\D_c$ is defined as
\begin{equation}
\lb{col}
(\D_c f)(v)=\sum\limits_{(v,\,u)\in\cA}\big(f(v)-f(u)\big), \quad f\in \ell^2(V), \quad v\in V.
\end{equation}
The combinatorial Schr\"odinger operator $\D_c+Q$ with a periodic potential $Q$ on periodic graphs was studied in \cite{KS14}. We note that in the case of a graph
with all vertices having the same degree $\vk_+$, the normalized Laplacian $\D$ and the combinatorial
Laplacian $\D_c$ (and, consequently, their spectra) are related by the simple identity
$\D_c=\vk_+\D$. However, in the case of an arbitrary graph the spectra of these
operators, in spite of many similar properties, may have significant differences.
For instance, the absolutely continuous spectrum of the normalized Laplacian $\D$
on the simplest periodic graph obtained from the square lattice by adding $N$ vertices
on each its edge has the form
$$
\s_{ac}(\D)=[0,2]
$$
(see Proposition \ref{TG2}). In the case of the combinatorial
Laplacian $\D_c$ on the same graph with $N=2$ the absolutely continuous spectrum
consists of three spectral bands separated by gaps (see \cite[p.~600]{KS14}).

\medskip

Results about Laplacians on periodic graphs are used in
spectral analysis of the Schr\"odinger operator with a decaying potential and
also for the study of the Laplacians on periodic graphs with various defects.
We briefly describe these works. The scattering problem for the Schr\"odinger
operator with a decaying potential on the lattice $\Z^d$,
$d>1$, was considered in the papers \cite{BS99,IK12,IM14,Ko10,KM17,RS09,SV01} (see also
references therein). The scattering on other graphs was studied in
\cite{A12,KS15,KMR18,PR18}. The Schr\"odinger operator
with a potential periodic in some directions
and finitely supported in other directions on arbitrary periodic graphs was investigated
in \cite{KS17a}. In \cite{AIM16,Ku14,KS17b,SS17} the Laplace
and Schr\"odinger operators on periodic graphs with different defects were considered.

\subsection{Edge indices}
In order to formulate our results we need to define the notion of the {\it edge index},
which was introduced in \cite{KS14}. Indices are important to study the
spectrum of Laplace and Schr\"odinger operators on periodic graphs, since
fiber operators are expressed in terms of edge indices of the fundamental graph
$\G_*=(V_*,\cE_*)$ (see the formula \eqref{l2.13}).

Let $\nu=\#V_*$, where $\#A$
is the number of elements in a set $A$. We fix any $\nu$ vertices of the periodic
graph $\G$, which are not $\Z^d$-equivalent to one another and denote this vertex set
by $V_0$. We will call the set $V_0$ a \emph{fundamental vertex set of the graph $\G$.}
The set $V_0$ can be chosen by different ways. However, it is natural to choose
this set in the following way. By Lemma \ref{Lst}.i, there exists a subgraph
$T=(V_T,\cE_T)$ of the periodic graph $\G$ satisfying the following conditions:

1) \emph{$T$ is a tree, i.e., a connected graph without cycles{\rm;}}

2) {\it the set $V_T$ consists of $\n$ vertices of $\G$ that are not
$\Z^d$-equivalent to each other.}\\
From now on we assume that the fundamental vertex set $V_0$
coincides with the vertex set $V_T$ of the tree $T$. We note that the graph
$T$ is not unique (see Lemma \ref{Lst}.i).

For any vertex $v\in V$ the following unique representation holds true:
\begin{equation}
\lb{Dv} v=v_0+[v], \quad v_0\in
V_0,\quad [v]\in\Z^d.
\end{equation}
In other words, each vertex $v$ of the periodic graph $\G$ can be obtained from
a vertex $v_0\in V_0$ by the shift by an integer vector $[v]\in \Z^d$.
We will call the vector $[v]$ the \emph{coordinates of the vertex $v$ with
respect to the fundamental vertex set $V_0$}. For any oriented edge $\be=(u,v)\in\cA$
we define the {\bf edge index} $\t(\be)$ as the integer vector given by
\begin{equation}
\lb{in}
\t(\be)=[v]-[u]\in\Z^d,
\end{equation}
where, due to \eqref{Dv}, we have
$$
u=u_0+[u],\quad v=v_0+[v], \quad u_0,v_0\in V_0,\quad [u],[v]\in\Z^d.
$$
For example, for the graph $\G$ shown in Fig.\ref{ff.0.11} the index of the edge
$(v_3,v_1+a_1+a_2)$ is equal to $(1,1)$ and the edge $(v_1,v_4)$ has a zero index.
Generally speaking, edge indices depend on the choice of the set $V_0$.

We define a surjection $\gf_\cA:\cA\rightarrow\cA_*=\cA/\Z^d$, which maps each oriented
edge $\be\in\cA$ of $\G$ to its equivalence class, i.e., an oriented edge
$\be_*=\gf_{\cA}(\be)$ of the fundamental graph $\G_*$. For each edge $\be_*\in\cA_*$
we define the edge index $\t(\bf e_*)$ in the following way:
\begin{equation}
\lb{inf}
\t(\bf e_*)=\t(\be).
\end{equation}

\setlength{\unitlength}{1.0mm}
\begin{figure}[h]
\centering
\unitlength 1mm 
\linethickness{0.4pt}
\ifx\plotpoint\undefined\newsavebox{\plotpoint}\fi 

\begin{picture}(60,50)(0,0)

\multiput(10,10)(4,0){10}{\line(1,0){2}}
\multiput(10,30)(4,0){10}{\line(1,0){2}}
\multiput(10,50)(4,0){10}{\line(1,0){2}}

\multiput(10,10)(0,4){10}{\line(0,1){2}}
\multiput(30,10)(0,4){10}{\line(0,1){2}}
\multiput(50,10)(0,4){10}{\line(0,1){2}}

\put(10,10){\vector(1,0){20.00}}
\put(10,10){\vector(0,1){20.00}}

\put(20,8){$\scriptstyle a_1$}
\put(6,20){$\scriptstyle a_2$}

\put(15,15){\line(1,1){20.00}}
\put(15,15.2){\line(1,1){10.00}}
\put(15,15.4){\line(1,1){10.00}}
\put(15,14.9){\line(1,1){10.00}}
\put(15,14.8){\line(1,1){10.00}}

\put(25,15){\line(-1,1){10.00}}
\put(25,15.2){\line(-1,1){10.00}}
\put(25,15.4){\line(-1,1){10.00}}
\put(25,14.9){\line(-1,1){10.00}}
\put(25,14.8){\line(-1,1){10.00}}

\put(15,15){\line(1,0){20.00}}
\put(15,15){\line(0,1){20.00}}

\put(15,15){\circle*{1}}
\put(15,25){\circle*{1}}
\put(20,20){\circle*{1}}
\put(25,25){\circle*{1}}
\put(25,15){\circle*{1}}
\put(12,13){$\scriptstyle v_1$}
\put(12,26.5){$\scriptstyle v_2$}\put(33,13){$\scriptstyle v_1+a_1$}
\put(24,12.5){$\scriptstyle v_4$}\put(35,32.5){$\scriptstyle v_1+a_1+a_2$}
\put(23,26.5){$\scriptstyle v_3$}\put(12,36){$\scriptstyle v_1+a_2$}
\put(18.5,22){$\scriptstyle v_5$}


\put(35,15){\circle{0.7}}
\put(35,25){\circle{0.7}}
\put(40,20){\circle{0.7}}
\put(45,25){\circle{0.7}}
\put(45,15){\circle{0.7}}

\put(15,35){\circle{0.7}}
\put(15,45){\circle{0.7}}
\put(20,40){\circle{0.7}}
\put(25,45){\circle{0.7}}
\put(25,35){\circle{0.7}}


\put(35,35){\circle{0.7}}
\put(35,45){\circle{0.7}}
\put(40,40){\circle{0.7}}
\put(45,45){\circle{0.7}}
\put(45,35){\circle{0.7}}
\end{picture}

\vspace{-0.5cm} \caption{ \footnotesize A graph $\G$ with the fundamental vertex set
$\{v_1,\ldots,v_5\}$; only edges of the fundamental graph $\G_*$ are shown in the
figure; the vectors $a_1,a_2$ produce an action of the group $\Z^2$;
edges of the tree $T$ are marked with bold lines.} \label{ff.0.11}
\end{figure}

In other words, edge indices of the fundamental graph $\G_*$ are induced by
indices of the corresponding edges of the periodic graph $\G$.
The index of a fundamental graph edge with respect to the fixed fundamental vertex set
$V_0$ is uniquely determined by formula \eqref{inf}, since
$$
\t(\be+m)=\t(\be),\quad \forall\, (\be,m)\in\cA\ts\Z^d.
$$
From the definition \eqref{in} of the edge index it follows that
\textbf{all edges of the tree $T$ have zero indices}, i.e.,
\begin{equation}
\lb{iesp}
\t(\be)=0, \quad \forall \ \be\in\cE_{T}.
\end{equation}
Edges with non-zero indices exist on any periodic graph and provide its
connectivity. We denote by $\cB$ and $\cB_*$ the sets of all edges with non-zero
indices of the periodic graph $\G$ and the fundamental graph $\G_*$, respectively.

\subsection
{The direct integral and the spectrum of the Schr\"odinger operator}
It is well known that periodic operators can be decomposed into a direct integral
(for the continuous case see \cite{RS78}). The existence of the direct integral
\eqref{raz} for the discrete Schr\"odinger operator on periodic graphs was discussed
in many papers, see \cite{KSS98,Ku89,RR07}. In the case of a concrete periodic graph
it is not difficult to write down an explicit expression for the fiber operator
(see, e.g., Subsection 4.2.2 in \cite{BK13} and also \cite{HKSW07}). In the case of
an arbitrary periodic graph explicit forms of the fiber operators each of which acts
in the same space $\ell^2(V_*)$, are given in \cite{KS14,KS16a,S13}.
Repeating the proof of Theorem 1.1.i for the combinatorial Laplacian in \cite{KS14},
we obtain the following statement.

\begin{proposition}\lb{diso}
The Schr\"odinger operator $H=\D+Q$ in the space $\ell^2(V)$ is decomposed into
a constant fiber direct integral
\begin{equation}
\lb{raz}
\mH={1\over(2\pi)^d}\int\limits_{\T^d}^{\os}\ell^2(V_*)\,d\vt, \quad  UH U^{-1}={1\over(2\pi)^d}\int\limits^\oplus_{\T^d}H(\vt)d\vt,
\end{equation}
where $\T^d=\R^d/(2\pi\Z)^d$, $U:\ell^2(V)\to\mH$ is some unitary operator
(the Gelfand transformation). Here the fiber Schr\"odinger operator $H(\vt)$ and the
fiber Laplacian $\D(\vt)$ for each $\vt\in \T^d$ have the form
\begin{equation}
\label{Hvt}
H(\vt)=\D(\vt)+Q,
\end{equation}
\begin{equation}
\label{l2.13}
\big(\D(\vt)f\big)(v)=f(v)
-\sum_{\be=(v,\,u)\in\cA_*}\frac{e^{i\lan\t(\be),\,\vt\rangle}}{\sqrt{\vk_v\vk_u}}f(u), \quad v\in V_*,
\end{equation}
where $\t(\be)$ is the index of the edge $\be$, defined by the formulas \eqref{in},
\eqref{inf}{\rm;} $\lan\,\cdot\,,\,\cdot\,\rangle$ is the standard inner product in
$\R^d$, and $\vk_v$ is the degree of the vertex $v$.
\end{proposition}

\no\textbf{Remarks.} 1) The explicit form \eqref{l2.13} of the fiber operator $\D(\vt)$
is important to study spectral properties of Laplace and Schr\"odinger operators
on periodic graphs (see the proofs of Theorems {\ref{TMR1}--\ref{Tfb})}.

2) In \cite{S13} the fiber operator $\D(\vt)$ is expressed in terms of
coordinates of edges (considered as directed line segments connecting
the corresponding vertices) of the initial periodic graph realized in the space
$\R^d$. The coordinate vector of each edge

$\bu$ is a vector in the space $\R^d$ (not necessarily integer);

$\bu$ is a \textbf{non-zero} vector (except for loop edges of the periodic graph).\\
In contrast to the coordinate vector introduced in \cite{S13}, the edge index has
the following important properties:

$\bu$ the edge index is an \textbf{integer} vector in the space $\R^d$;

$\bu$ on the fundamental graph $\G_*=(V_*,\cE_*)$ there exist at least $\#V_*-1$ edges
with \textbf{non-zero indices} (see the formula \eqref{iesp}).

It should be noted that the number of the fundamental graph edges can be
arbitrarily large integer number, while the number of edges with non-zero indices
can be $d$. For such graphs the dependence of the fiber operator $\D(\vt)$ on the
quasimomentum $\vt$ in Sunada's decomposition is rather complicated, while the same
dependence in the expression \eqref{l2.13} is essentially simpler. We also note that
in the paper \cite{KS18} we obtained a decomposition of the Laplacian $\D$ on
periodic graphs into a direct integral where the fiber operator $\D(\vt)$ depends
on the minimal number of non-zero indices and showed that this number is an
invariant of the periodic graph.

\medskip

The decomposition \eqref{raz} and the standard theory of periodic operators (see
\cite{RS78}, Theorem XIII.85) describe the spectrum of the Schr\"odinger operator
$H=\D+Q$. Each fiber operator $H(\vt)$, $\vt\in\T^d$, has $\n=\# V_*$ eigenvalues
$\l_{n}(\vt)$, $n\in\N_\n=\{1,\ldots,\n\}$, labeled in non-decreasing order counting
multiplicity:
\begin{equation}
\label{eq.3} \l_{1}(\vt)\leq\l_{2}(\vt)\leq\ldots\leq\l_{\nu}(\vt),\quad \n=\# V_*,
\quad \forall\,\vt\in\T^d.
\end{equation}
Since the operator $H(\vt)$ is self-adjoint and analytic in $\vt\in\T^d$, each
$\l_{n}(\,\cdot\,)$, $n\in\N_\n$, is a real-valued and piecewise analytic function
on the torus $\T^d$ and defines the \emph{spectral band} $\s_n(H)$:
\begin{equation}
\lb{ban.1}
\s_n(H):=\s_{n}=[\l_{n}^-,\l_{n}^+]=\l_{n}(\T^d).
\end{equation}
Then the spectrum of the Schr\"odinger operator $H$ on the graph $\G$ is given by
\[\lb{spec}
\s(H)=\bigcup_{\vt\in\T^d}\s\big(H(\vt)\big)=\bigcup_{n=1}^{\nu}\s_n(H).
\]
Since the fiber operator $H(\vt)$ acts in a finite-dimensional space and is analytic in
$\vt\in\T^d$ (moreover, it is an entire function of $\vt\in \C^d$), according to the
general theory, see \cite{GN98}, the singular continuous spectrum of the operator $H$ is
absent. Note that if $\l_{n}(\,\cdot\,)=C_{n}=\const$ on some subset of $\T^d$ of
positive Lebesgue measure, then the operator~$H$ on the graph~$\G$ has the
eigenvalue~$C_{n}$ of infinite multiplicity. We call $\{C_{n}\}$ a \emph{flat band}.
Thus, the spectrum of the Schr\"odinger operator $H$ on the periodic graph $\G$ has
the form
\begin{equation}
\lb{r0}
\s(H)=\s_{ac}(H)\cup \s_{fb}(H).
\end{equation}
Here $\s_{ac}(H)$ is the absolutely continuous spectrum, which is a union of
non-degenerate bands from \eqref{spec}, and $\s_{fb}(H)$ is the set of all eigenvalues
of infinite multiplicity (flat bands). An open interval between two neighboring
non-degenerate spectral bands is called a \emph{gap}.

The eigenvalues of the fiber Laplacian $\D(\vt)$ will be denoted by $\l^0_{n}(\vt)$,
$n\in\N_\n$. The spectral bands $\s_n(\D)$, $n\in\N_{\n}$, for the Laplacian $\D$
have the form
\begin{equation}
\lb{ban0} \s_n(\D):=\s_{n}^0=[\l_{n}^{0-},\l_{n}^{0+}]=\l_{n}^0(\T^d).
\end{equation}

\section{Main results}\lb{Sec2}
\setcounter{equation}{0}

\subsection{Estimates for the first spectral band of the Schr\"odinger operator}
It is known (see \cite[Theorem~1]{SS92}) that
\begin{equation}
\lb{bes}
\l_1^-=\min\limits_{\vt\in\T^d}\l_1(\vt)=\l_1(0).
\end{equation}
The first band function $\l_1(\vt)$ has a Taylor series expansion about the point 0:
\begin{equation}\lb{tem}
\begin{aligned}
&\l_1(\vt)=\l_1(0)+\frac{1}{2}\sum_{i,j=1}^{d} M_{ij}\,\vt_i \vt_j+O(|\vt|^3),
\\
&\text{where}\quad M_{ij}={\pa^2\l_1(0)\over\pa\vt_i\pa\vt_j},\quad \vt=(\vt_j)_{j=1}^d.
\end{aligned}
\end{equation}
The entries of the matrix $m=M^{-1}$, where $M=\{M_{ij}\}$, represent a tensor, which
is called \emph{the effective mass tensor} at the bottom of the spectrum.
The effective mass approximation is a standard approach in solid state physics.
By this approach, for energies close to $\l_1(0)$ a complicated Hamiltonian is replaced
by the model Hamiltonian $-{\D\over 2m}$, where $\D$ is the Laplacian and $m$ is
the effective mass. In a general case the effective mass depends on the direction
in crystals  and represents a tensor.

In the paper \cite{KS16a} the authors obtained upper estimates for the effective masses
associated with the ends of each spectral band of the discrete Laplacians
on periodic graphs in terms of geometric parameters of the graphs. Moreover, in the
case of the bottom of the spectrum they determined two-sided
estimates on the effective mass in terms of geometric parameters of
the graphs. In the paper \cite{K08} the effective masses for magnetic Schr\"odinger operators on zigzag
nanotubes were estimated. We note that in the case of the Schr\"odinger operators
with periodic potentials in the space $\R^d$ the effective mass tensor was studied in
\cite{BS04,KiS87,Sh06}.

We formulate two-sided estimates on the lengths of the first spectral
bands and on the effective masses at the bottom of the spectrum of the normalized
Laplace and Schr\"odinger operators.

\begin{theorem}
\lb{TMR1}
i) The smallest eigenvalue $\l_1(0)$ of the operator $H(0)$ is simple,
and all components of a corresponding eigenvector $\p\in\ell^2(V_*)$
are strictly positive.

ii) The first band $\s_1(H)=[\l_1^-,\l_1^+]$ of the Schr\"odinger operator
$H=\D+Q$ is non-degenerate, i.e., $\l_1^-<\l_1^+$, and the following estimate
holds:
\[
\lb{eq.TMR1}
\begin{aligned}
&c_0^{-2}|\s_1(\D)|\leq|\s_1(H)|\leq
c_0^{2}\,|\s_1(\D)|,\\
&c_0={\psi_+\over \psi_-},\quad \psi_-=\min\limits_{v\in V_*}\frac{\psi(v)}{\sqrt{\vk_v}},\quad \psi_+=\max\limits_{v\in V_*}\frac{\psi(v)}{\sqrt{\vk_v}},
\end{aligned}
\]
where $\s_1(\D)$ is the first band of the Laplacian~$\D$, and $\vk_v$ is the degree
of the vertex~$v$.

iii) The effective mass tensors $m_0$ and $m$ at the bottom of the spectrum
of the Laplacian $\D$ and the Schr\"odinger operator $H=\D+Q$, respectively,
satisfy:
\begin{equation}
\lb{eem}
c_0^{-2}m_0\leq m\leq
c_0^{2}\, m_0.
\end{equation}
\end{theorem}

\no \textbf{Remarks.} 1) Theorem \ref{TMR1} remains true for the Schr\"odinger operator
$\D_{c}+Q$, where $\D_c$ is the combinatorial Laplacian defined by
\eqref{col}. In this case
$$
\psi_-=\min\limits_{v\in V_*}\psi(v),\quad \psi_+=\max\limits_{v\in V_*}\psi(v)
$$
and the proof repeats the proof of Theorem \ref{TMR1}.

2) For the lattice $\Z^d$, Theorem \ref{TMR1} was proved in \cite{KiS87}.
For an arbitrary periodic graph the estimates \eqref{eq.TMR1}, \eqref{eem} are
new.

\subsection{Estimate of the Lebesgue measure of the spectrum}
We estimate the Lebesgue measure of the spectrum of the Schr\"odinger operator
$H=\D +Q$ in terms of geometric parameters of the graph (the number of edges with
non-zero indices and vertex degrees) and the sum of gap lengths in terms of geometric
parameters of the graph and the potential $Q$.

\begin{theorem}
\lb{T2.2}
i) The Lebesgue measure $|\s(H)|$ of the spectrum of the Schr\"odinger operator
$H=\D+Q$ satisfies the inequality:
\begin{equation}
\lb{eq.7}
|\s(H)|\le\sum_{n=1}^{\n}|\s_n(H)|\le2\gz,\quad \gz=\sum_{v\in V_*}\frac{\gz_v}{\vk_v},
\end{equation}
where $\gz_v$ is the number of the fundamental graph edges having non-zero indices and
starting at the vertex $v\in V_*$, $\vk_v$ is the degree of the vertex $v$, and $\n=\#V_*$.
Moreover, if there exist $s$ gaps $\g_1(H),\ldots,\g_s(H)$ in the spectrum of
the operator $H$, then the following estimate holds true{\rm:}
\begin{equation}
\lb{GEga}
\begin{aligned}
&\sum_{n=1}^s|\g_n(H)|\ge
\l^+_{\n}-\l_{1}^--2\gz\ge C_0-2\gz,
\\
&C_0=|\l^{0+}_{\n}-q_\bu|, \quad q_\bu=\max_{v\in V_*} Q(v)-\min_{v\in V_*} Q(v),
\end{aligned}
\end{equation}
where $\l^{0+}_{\n}$ is the upper endpoint of the spectrum of the Laplacian $\D;$
and $\l^{-}_{1}$, $\l^{+}_{\n}$ are the lower and upper endpoints of the spectrum of
the Schr\"odinger operator $H$.

ii) The estimates \eqref{eq.7} and the first estimate in \eqref{GEga} become
identities for some classes of graphs, see \eqref{loo}.
\end{theorem}

\no \textbf{Remarks.} 1) We recall that the spectrum of the normalized Laplacian
$\s(\D)\ss[0,2]$. Therefore, the estimates \eqref{eq.7} for the Laplacian $\D$ are
effective if and only if $\gz<1$. This condition holds true when, for each vertex
$v\in V_*$, the number of edges having non-zero indices and starting at $v$ is
sufficiently small compared to the degree of the vertex.

2) The number $\gz$ satisfies the following simple estimates (see Lemma \ref{Lst}.ii):
$$
\gz\leq 1, \quad \textrm{if} \quad \n=1;\quad
\gz\leq \n-\sum_{v\in V_*}{1\over \vk_v}, \quad \textrm{if} \quad \n\geq2;
\quad \n=\#V_*.
$$

\medskip

\subsection{The Schr\"odinger operator on loop graphs}
A periodic graph~$\G$ is called a {\it loop  graph} if each edge
of the fundamental graph
$$
\G_*=(V_*,\cE_*)
$$
with a non-zero index is a loop, i.e., the edge has the form $(v,v)$ for some vertex
$v\in V_*$.

A loop graph $\G$ is called a {\it precise} loop graph if
$$
\cos\lan\t(\be),\vt_0\rangle=-1
$$
for all edges $\be\in\cB_*$ with non-zero indices and some $\vt_0\in\T^d$, where
$\t(\be)$ is the index of the edge $\be$. The point $\vt_0$ is called a {\it precise
quasimomentum} of the loop graph $\G$.

\medskip

\no \textbf{Remark.} The class of loop graphs is large enough. The simplest
example of a precise loop graph is the $d$-dimensional lattice. More
complicated examples of loop graphs are discussed in Proposition 2.3 in \cite{KS14}.

\medskip

We describe all spectral bands for the Schr\"odinger operator on precise loop graphs.

\begin{theorem}\lb{T100}
i) Let $\G$ be a loop graph. Then the lower endpoints of the spectral bands
$\s_n=\s_n(H)=[\l_n^-,\l_n^+]$ of the Schr\"odinger operator $H$ on the graph $\G$
satisfy
\begin{equation}
\lb{eq.5} \l_n^-=\l_n(0),\quad \forall \;n\in \N_\n.
\end{equation}

ii) Let $\G$ be a precise loop graph with a precise quasimomentum
$\vt_0\in\T^d$. Then
\begin{equation}
\lb{es}
\s_n=[\l_n^-,\l_n^+]=[\l_n(0),\l_n(\vt_0)],\quad \forall\,n\in \N_\n,
\end{equation}
\begin{equation}
\lb{es.1}
\sum_{n=1}^\n|\s_n|=2\gz,
\end{equation}
where $\gz$ is defined in \eqref{eq.7}. In particular, if all edges of the
fundamental graph $\G_*=(V_*,\cE_*)$ with non-zero indices have the form $(v,v)$
for some vertex $v\in V_*$, then
\begin{equation}
\lb{loo}
|\s(H)|=\sum_{n=1}^\n|\s_n|=2\gz.
\end{equation}
\end{theorem}

\no \textbf{Remarks.} 1) By \eqref{es.1}, the total length of the spectral bands
of the Schr\"odinger operator $H=\D+Q$ on precise loop graphs
does not depend on the potential $Q$.

2) $\l_n^-$, $n\in\N_\n$, are the eigenvalues of the Schr\"odinger operator $H(0)$
defined by the formulas \eqref{Hvt}, \eqref{l2.13}, on the fundamental graph $\G_*$.
Identities similar to \eqref{es} hold in the case of $N$-periodic Jacobi matrices
on the lattice $\Z$ (and for the Hill operator). The spectrum of these operators is
absolutely continuous and is a union of spectral bands separated by gaps. The endpoints
of the bands are the so-called $2N$-periodic eigenvalues.

\medskip

\subsection{Spectrum of the Laplacian on specific periodic graphs}
We discuss the possible number of non-degenerate spectral bands and the possible number
of eigenvalues of infinite multiplicity (flat bands) for the normalized Laplacian.
According to the standard methods of analytic continuation (see, e.g.,
\cite{RS78,T73}), $\l_*$ is an eigenvalue of the operator $H$ if and only if
$\l_*=\const$ is an eigenvalue of the operator $H(\vt)$ for all $\vt\in\T^d$.
We define the multiplicity of a flat band in the following way:
a flat band $\{\l_*\}$ of the operator $H$ has the multiplicity $m$ if and only if
$\l_*=\const$ is an eigenvalue of the operator $H(\vt)$ of multiplicity $m$ for
almost all $\vt\in\T^d$.

\begin{theorem}\lb{Tfb}
Let $\n,d\geq2$. Then the following statements hold true.

i) There exists a $\Z^d$-periodic graph $\G$ such that the spectrum
of the Laplacian $\D$ on $\G$ consists of two non-degenerate spectral bands
and $\n-2$ flat bands (counting multiplicity) lying in the gap.

ii) There exists a $\Z^d$-periodic graph $\G$ such that the spectrum
of the operator $\D$ on $\G$ has $[\frac{\n-1}d]$ distinct flat bands, each of which
has multiplicity $d-1$ and is embedded into the absolutely continuous spectrum
$\s_{ac}(\D)=[0,2]$.

iii) There exists a $\Z^d$-periodic graph $\G$ such that the Lebesgue measure
of the spectrum of the Schr\"odinger operator $H=\D+Q$ on $\G$ satisfies
$|\s(H)|={4d\over \n+2d-1}$. In particular, we have $|\s(H)|\rightarrow0$ as $\n\to\infty$.
\end{theorem}

\no \textbf{Remark.} The eigenspace corresponding to an eigenvalue of infinite multiplicity
is generated by finitely supported eigenfunctions (see, e.g.,
\cite[Theorem 4.5.2]{BK13}). We note that this statement holds not only in the case
of periodic graphs, but also for graphs with a more general structure, see \cite{V05}.

\medskip

We briefly describe the layout of the paper. Auxiliary Section \ref{Sec3} is devoted to
a factorization of the fiber Laplacian. In Section \ref{Sec4} we prove Theorems \ref{TMR1}
and \ref{T2.2} about spectral estimates for the Schr\"odinger operator $H$. In
Section \ref{Sec5} we describe the spectrum of the operator $H$ on loop graphs and prove
Theorem \ref{T100}. In that section we also study spectral properties of the Laplacian
$\D$ on bipartite graphs. In Sections \ref{Sec6}, \ref{Sec7} the spectrum of the Laplace and
Schr\"odinger operators on some perturbations of the $d$-dimensional lattice by adding
vertices and edges (in a periodic way) is described. Theorem \ref{Tfb} is also
proved there. In Section \ref{Sec8} we recall some well-known properties of matrices needed
to prove our main results.

\section{\lb{Sec3} Factorization of the Laplace operator}
\setcounter{equation}{0}

\subsection{The existence of the tree $T$ and estimates for the number $\gz$}

\begin{lemma}\label{Lst}
i) There exists a subgraph $T=(V_T,\cE_T)$ of the periodic graph~$\G$,
satisfying the following conditions:

1) $T$ is a tree, i.e., a connected graph without cycles;

2) the set $V_T$ consists of $\n$ vertices of the graph $\G$, which are not
$\Z^d$-equivalent to each other.

ii) For the number $\gz$ defined in \eqref{eq.7}, the following estimates
hold:
$$
\gz\leq 1, \quad \textrm{if} \quad \n=1;\quad
\gz\leq \n-\sum_{v\in V_*}{1\over \vk_v}, \quad \textrm{if} \quad \n\geq2;
\quad \n=\#V_*,
$$
where $\vk_v$ is the degree of the vertex $v$.
\end{lemma}

\no \textbf{Remark.} Due to the construction, the graph $T$ is not unique.

\medskip

\no \textbf{Proof.} \emph{i}) We describe a construction of the graph
$T=(V_T,\cE_T)$. First, we put $V_T=\varnothing$. Next, let $v$ be an arbitrary
vertex of the graph~$\G$. We add it to the set $V_T:=V_T\cup\{v\}$. Since
$\G$ is connected, there exist vertices adjacent to~$v$. We add to the set~$V_T$
those of them which are not $\Z^d$-equivalent to the vertices from~$V_T$.
We proceed similarly with each vertex newly added to~$V_T$ and so on. Since the number
of the fundamental graph vertices is finite, this process of adding vertices will
finish after a finite number of steps. The obtained set~$V_T$ satisfies
the condition 2). Indeed, assume that some vertex $u\in V\sm V_T$ is not
$\Z^d$-equivalent to any vertex from~$V_T$. Due to the construction
of the set~$V_T$ and the periodicity of the graph $\G$, this means that the vertex~$u$
is adjacent neither to vertices from the set~$V_T$, nor to their equivalent vertices.
This contradicts the connectivity of the periodic graph. Thus, the condition 2) holds.

Let $\G_1=(V_T,\cE_1)$ be a subgraph of $\G$ with the vertex set $V_T$ and
the edge set $\cE_1$, which contains only edges of $\G$ with both endpoints in $V_T$.
Due to the construction of the set $V_T$, the graph $\G_1$ is connected. Then the
spanning tree of the graph $\G_1$, i.e., a tree containing all vertices of the graph
$\G_1$, is a required subgraph of $T$.

\emph{ii}) Let $\n=1$, i.e., $V_*$ consists of a single vertex $v$. Then, by
\eqref{eq.7}, $\gz=\frac{\gz_v}{\vk_v}\leq 1$.

Now let $\n\geq2$. Due to the construction of the graph $T$, for each vertex $v$ in the
fundamental set $V_0=V_T$, there exists at least one edge of the tree $T=(V_T,\cE_T)$
incident to $v$. Then, by \eqref{inf} and \eqref{iesp}, for each vertex $v\in V_0$ and,
consequently, for each vertex $v$ of the fundamental graph, the number $\gz_v$ of
edges having non-zero indices and starting at $v$ is not greater than $\vk_v-1$.
This yields the required estimate:
$$
\gz=\sum_{v\in V_*}\frac{\gz_v}{\vk_v}\leq\sum_{v\in V_*}\frac{\vk_v-1}{\vk_v}=
\n-\sum_{v\in V_*}{1\over \vk_v}\,.\qqq \BBox
$$

\subsection{Factorization of the fiber Laplacian} We introduce the Hilbert space
\begin{equation}
\lb{el2A}
\ell^2(\cA)=\{\phi\colon \cA\ra\C\mid \phi(\ul\be)=-\phi(\be)
\textrm{ for all } \be\in\cA \; \textrm{ and } \, \lan
\phi,\phi\rangle_{\cA}<\iy\},
\end{equation}
with the inner product
\begin{equation}
\lb{incA}
\lan
\phi_1,\phi_2\rangle_{\cA}={1\over 2}\sum_{\be\in\cA}\phi_1(\be)\ol{\phi_2(\be)}.
\end{equation}

It is known (see, e.g., \cite{MW89}), that the Laplacian $\D$ has the following
factorization:
\begin{equation}
\lb{falo}
\D=\na^*\na,
\end{equation}
where the operator $\na:\ell^2(V)\ra \ell^2(\cA)$ is given by
\begin{equation}
(\na f)(\be)={f(v)\over \sqrt{\vk_v}}-{f(u)\over \sqrt{\vk_u}},\quad
f\in \ell^2(V),\quad \be=(v,u)\in\cA.
\end{equation}
The conjugate operator $\na^*\colon\ell^2(\cA)\ra \ell^2(V)$ has the form
$$
(\na^*
\phi)(v)={1\over \sqrt{\vk_v}}\sum_{\be=(v,\,u)\in\cA}\phi(\be),\quad \phi\in \ell^2(\cA),\quad v\in V.
$$
The quadratic form $\lan\D f,f\rangle_V$ of the Laplacian $\D$ is given by
\begin{equation}
\lb{qflo1}
\lan \D f,f\rangle_V={1\over 2}\sum_{(v,\,u)\in\cA}
\bigg|{f(v)\over \sqrt{\vk_v}}-{f(u)\over \sqrt{\vk_u}}\,\bigg|^2,
\end{equation}
where $\lan\,\cdot\,,\,\cdot\,\rangle_V$ is the inner product in the space $\ell^2(V)$:
\begin{equation}
\lb{ipV}
\lan
f_1,f_2\rangle_{V}=\sum_{v\in V}f_1(v)\ol{f_2(v)}, \quad \forall\, f_1,f_2\in\ell^2(V).
\end{equation}

We obtain a similar representation for the fiber Laplacian~$\D(\vt)$.

\begin{theorem}\label{FOFa}
i) For each $\vt\in\T^d$ the fiber Laplacian $\D(\vt)$ defined by the formula
\eqref{l2.13} satisfies the identity
\begin{equation}
\lb{fact}
\D(\vt)=\na^*(\vt)\na(\vt),
\end{equation}
where the operator $\na(\vt):\ell^2(V_*)\ra \ell^2(\cA_*)$ is given by
\begin{equation}
\begin{aligned}
\lb{navt}
\big(\na(\vt)
f\big)(\be)&=e^{-i\lan\t(\be),\vt\rangle/2}\,{f(v)\over \sqrt{\vk_v}}-
e^{i\lan\t(\be),\vt\rangle/2}\,{f(u)\over \sqrt{\vk_u}},
\\
 f&\in \ell^2(V_*), \quad \be=(v,u)\in\cA_*,
\end{aligned}
\end{equation}
$\t(\be)$ is the index of the edge $\be$ defined by the formulas \eqref{in},
\eqref{inf}. The conjugate operator $\na^*(\vt)\colon\ell^2(\cA_*)\ra \ell^2(V_*)$
has the form
\begin{equation}
\lb{coop}
(\na^*(\vt)
\phi)(v)\!=\!{1\over \sqrt{\vk_v}}\!\!\sum_{\be=(v,\,u)\in\cA_*}\!\!
e^{i\lan\t(\be),\vt\rangle/2}\phi(\be),\quad \phi\in\ell^2(\cA_*),\quad v\in V_*.
\end{equation}

ii) The quadratic form $\lan\D(\vt) f,f\rangle_{V_*}$ of the fiber Laplacian
satisfies the identity
\begin{equation}
\lb{qflo}
\lan \D(\vt)
f,f\rangle_{V_*}={1\over 2}\sum_{\be=(v,\,u)\in\cA_*}
\bigg|{f(v)\over \sqrt{\vk_v}}-
e^{i\lan\t(\be),\vt\rangle}\,{f(u)\over \sqrt{\vk_u}}\,\bigg|^2,
\end{equation}
where summation is over all oriented edges $\be\in\cA_*$.
\end{theorem}

\no \textbf{Proof.} Let $\vt\in\T^d$, $f\in\ell^2(V_*)$,
$\phi\in\ell^2(\cA_*)$. We denote $\wt f(v)={f(v)\over \sqrt{\vk_v}}$ for all
$v\in V_*$.

\emph{i}) First, we prove \eqref{coop}. Using the identity $\t(\be)=-\t(\ul\be)$ for all
$\be\in\cA_*$ and the formulas \eqref{el2A}, \eqref{incA}, \eqref{navt}, we obtain
\begin{multline}\lb{coop1}
\lan\na(\vt)f,\phi\rangle_{\cA_*}=
{1\over 2}\sum_{\be\in\cA_*}\big(\na(\vt)f\big)(\be)\ol{\phi(\be)}
\\
={1\over 2}\sum_{\be=(v,\,u)\in\cA_*}\big(e^{-i\lan\t(\be),\vt\rangle/2}\wt f(v)-
e^{i\lan\t(\be),\vt\rangle/2}\wt f(u)\big)\ol{\phi(\be)}
\\
={1\over 2}\sum_{\be=(v,\,u)\in\cA_*}
e^{-i\lan\t(\be),\vt\rangle/2}\wt f(v)\ol{\phi(\be)}
+{1\over 2}\sum_{\ul\be=(u,\,v)\in\cA_*}
e^{-i\lan\t(\ul\be),\vt\rangle/2}\wt f(u)\ol{\phi(\ul\be)}
\\
=\sum_{\be=(v,\,u)\in\cA_*}
e^{-i\lan\t(\be),\vt\rangle/2}\wt f(v)\ol{\phi(\be)}.
\end{multline}
On the other hand, due to \eqref{ipV}, \eqref{coop}, we have
\begin{multline}
\lb{coop2}
\lan f,\na^*(\vt)\phi\rangle_{V_*}=\sum_{v\in
V_*}f(v)\ol{(\na^*(\vt)\phi)(v)}\\=\sum_{v\in
V_*}\wt f(v)\sum_{\be=(v,\,u)\in\cA_*}
e^{-i\lan\t(\be),\vt\rangle/2}\,\ol{\phi(\be)}
=\sum_{\be=(v,u)\in\cA_*}
e^{-i\lan\t(\be),\vt\rangle/2}\wt f(v)\ol{\phi(\be)}.
\end{multline}
Comparing \eqref{coop1} and \eqref{coop2}, we see that the operator $\na^*(\vt)$
defined by the formula \eqref{coop} is indeed conjugate of $\na(\vt)$.

Next, we prove the identity \eqref{fact}. Using \eqref{navt}, \eqref{coop} and
\eqref{l2.13}, we obtain
\begin{multline*}
\big(\na^*(\vt)\na(\vt)f\big)(v)=
\sum_{\be=(v,\,u)\in\cA_*}{1\over \sqrt{\vk_v}}\,e^{i\lan\t(\be),\vt\rangle/2}
\big(\na(\vt)f\big)(\be)\\
=\sum_{\be=(v,\,u)\in\cA_*}{1\over \sqrt{\vk_v}}\,
e^{i\lan\t(\be),\vt\rangle/2}
\big(e^{-i\lan\t(\be),\vt\rangle/2}\wt f(v)-e^{i\lan\t(\be),\vt\rangle/2}\wt f(u)\big)
\\=
\sum_{\be=(v,\,u)\in\cA_*}{1\over \sqrt{\vk_v}}\,
\big(\wt f(v)-e^{i\lan\t(\be),\vt\rangle}\wt f(u)\big)
\\=
f(v)-\sum_{\be=(v,\,u)\in\cA_*}\frac{e^{i\lan\t(\be),\,\vt\rangle}}{\sqrt{\vk_v\vk_u}}f(u)
=\big(\D(\vt)f\big)(v),\quad
\forall v\in V_*.
\end{multline*}

\emph{ii}) From the identities \eqref{fact}, \eqref{incA}, \eqref{navt} it follows that
\begin{multline*}
\lan \D(\vt) f,f\rangle_{V_*}=\lan
\na(\vt)f,\na(\vt)f\rangle_{\cA_*}=
{1\over 2}\sum_{\be\in\cA_*}
\big|\big(\na(\vt)f\big)(\be)\big|^2
\\
={1\over 2}\sum_{\be=(v,u)\in\cA_*}
\big|e^{-i\lan\t(\be),\vt\rangle/2}\wt f(v)-e^{i\lan\t(\be),\vt\rangle/2}\wt f(u)\big|^2
\\
={1\over 2}\sum_{\be=(v,\,u)\in\cA_*}
\big|\wt f(v)-e^{i\lan\t(\be),\vt\rangle}\wt f(u)\big|^2.\qqq \BBox
\end{multline*}

\section{\lb{Sec4} Proof of the main results}
\setcounter{equation}{0}

\subsection{Estimates for the first spectral band of the Schr\"odinger operator}
We introduce the standard orthonormal basis of the space $\ell^2(V_*)$:
\begin{equation}
\lb{sonb}
(\gf_u)_{u\in V_*}, \quad \textrm{ where } \quad \gf_u(v)=\d_{uv}, \quad v\in V_*,
\end{equation}
$\d_{uv}$ is the Kronecker delta.

\begin{lemma}\label{TCo2}
For each $\vt\in \T^d$ the fiber Laplacian defined by the formula \eqref{l2.13}
in the standard basis \eqref{sonb} has the form
$$
\D(\vt)=\{\D_{uv}(\vt)\}_{u,v\in V_*},
$$
where
\begin{equation}
\lb{Devt}
\D_{uv}(\vt)=\begin{cases}
1-{1\over \vk_u}\sum\limits_{\be=(u,u)\in\cA_*}\cos\lan\t(\be),\vt\rangle, & \textrm{ if } \qq u=v;
\\
{-1\over \sqrt{\vk_u\vk_v}}
\sum\limits_{\be=(u,v)\in\cA_*} e^{-i\lan\t(\be),\,\vt\rangle},  &\textrm{ if } \qq u\sim v,\qq u\neq v;\\
  0, \qq & \textrm{ otherwise.}
\end{cases}
\end{equation}

\medskip
\noindent
Here ${\vk}_v$ is the degree of the vertex $v$, $\t(\be)$ is the index of the edge
$\be$ defined by the formulas \eqref{in}, \eqref{inf}.
\end{lemma}

\no \textbf{Proof.} Substituting the formula \eqref{l2.13} into the identity
$$
\D_{uv}(\vt)=\lan \gf_u,\D(\vt)\gf_v\rangle_{V_*}
$$
and using the fact that for each loop $\be=(u,u)\in\cA_*$ with the index $\t(\be)$
there exists a loop $\ul\be=(u,u)\in\cA_*$ with the index $-\t(\be)$, and the identity
$$
e^{-i\t(\be)}+e^{i\t(\be)}=2\cos\t(\be),
$$
we obtain \eqref{Devt}. \qq $\BBox$

\medskip

\no \textbf{Proof of Theorem \ref{TMR1}.} \emph{i}) The proof of this item can be found
in the paper \cite{SS92}. For the reader's convenience we give these simple arguments.
By \eqref{Devt}, the Laplacian on the fundamental graph $\G_*$ in the standard basis
\eqref{sonb} has the form
$
\D(0)=\{\D_{uv}(0)\}_{u,v\in V_*},
$
where
\begin{equation}
\lb{Devt1}
\D_{uv}(0)=\d_{uv}-{\vk_{uv}\over \sqrt{\vk_u\vk_v}},
\end{equation}
$\d_{uv}$ is the Kronecker delta, $\vk_{uv}$ is the number of edges of the form $(u,v)$
on the fundamental graph $\G_*$, and ${\vk}_v$ is the degree of the vertex $v$. Let
$$
A=\1_\n-\D(0),
$$
where $\1_\n$ is the identity $(\n\ts\n)$-matrix, $\n=\# V_*$. By~\eqref{Devt1},
all entries of the matrix $A$ are non-negative. Since the fundamental graph $\G_*$
is connected, by the matrix property iv), the matrix $A$ is irreducible.
The Schr\"odinger operator on the fundamental graph $\G_*$ in the standard basis
\eqref{sonb} has the form
 $$
 H(0)=D-A,
 $$
where
  $$
  D=\diag\big(Q(v)+1\big)_{v\in V_*}
  $$
is a diagonal matrix. Then, for sufficiently large $t\in\R$,
$$
-H(0)+t\1_\n=A-D+t\1_\n
$$
is an irreducible matrix with non-negative entries. Applying the matrix property v) to this
matrix, we obtain that the largest eigenvalue of this matrix is simple and there exists
a corresponding eigenvector with positive components. This yields
the required statement. \qq $\BBox$

\medskip

To prove the remaining items of Theorem \ref{TMR1}, we need the following lemma.

\begin{lemma}
\lb{T1}
Let $\p\in\ell^2(V_*)$ be an eigenvector  with positive components corresponding to
the smallest eigenvalue $\l_1(0)$ of the operator $H(0)$. Then

i) for any function $f\in\ell^2(V_*)$ and all $\vt\in\T^d$ the following
identity holds
\begin{equation}
\lb{eq.L1}
\begin{aligned}
\big\lan\big(H(\vt)-\l_1(0)\1\big)\p f,\p f\big\rangle_{V_*}
={1\over 2}\sum_{\be=(v,\,u)\in\cA_*}
c_{uv}\big|f(v)&-e^{i\lan\t (\be),\,\vt\rangle}f(u)\big|^2,
\\
 c_{uv}&=\frac{\p(u)\p(v)}{\sqrt{\vk_u\vk_v}},
\end{aligned}
\end{equation}
where $\1$ is the identity operator, and $\t(\be)$ is the index of the edge $\be\in\cA_*$
defined by the formulas \eqref{in}, \eqref{inf};

ii) the following estimate holds:
\begin{equation}
\lb{eq.T1}
\begin{aligned}
&c_0^{-2}\l_1^{(0)}(\vt)\leq\l_1(\vt)-\l_1(0)\leq
c_0^2\,\l_1^{(0)}(\vt), \quad
\forall\,\vt\in\T^d,\\
&c_0={\psi_+\over \psi_-},\quad \psi_-=\min\limits_{v\in V_*}\frac{\psi(v)}{\sqrt{\vk_v}},\quad \psi_+=\max\limits_{v\in V_*}\frac{\psi(v)}{\sqrt{\vk_v}},
\end{aligned}
\end{equation}
where $\l_1(\vt)$ and $\l_1^{(0)}(\vt)$ are the smallest eigenvalues of the operators
$H(\vt)$ and $\D(\vt)$, respectively.
\end{lemma}

\no \textbf{Proof.} \emph{i}) We use some arguments from \cite{KiS87}. Let $f\in\ell^2(V_*)$.
From the identities \eqref{ipV} and $H(0)\p=\l_1(0)\p$ it follows that
$$
\lan\l_1(0)\p f,\p f\rangle_{V_*}=\sum_{v\in V_*}\l_1(0)\p^2(v)|f(v)|^2
=\sum_{v\in V_*}\big(H(0)\p\big)(v)\p(v)|f(v)|^2.
$$
Then for each $\vt\in\T^d$, using \eqref{Hvt}, we have
\begin{multline}
\lb{eq.L2}
\big\lan\big(H(\vt)-\l_1(0)\1\big)\p f,\p f\big\rangle_{V_*}
=\sum_{v\in V_*}(H(\vt)\p f)(v)\p(v)\ol f(v)-\sum_{v\in V_*}\big(H(0)\p\big)(v)\p(v)|f(v)|^2
\\=\sum_{v\in V_*}(\D(\vt)\p f)(v)\p(v)\ol f(v)-\sum_{v\in V_*}\big(\D(0)\p\big)(v)\p(v)|f(v)|^2.
\end{multline}
Substituting the expression \eqref{l2.13} for the operator $\D(\vt)$ into the last
identity, we obtain
\begin{multline}
\lb{eq.L22}
\big\lan\big(H(\vt)-\l_1(0)\1\big)\p f,\p f\big\rangle_{V_*}
=\sum_{v\in V_*}\Big(\p(v)f(v)-
\sum_{\be=(v,\,u)\in\cA_*}\frac{e^{i\lan\t(\be),\,\vt\rangle}}
{\sqrt{\vk_v\vk_u}}\p(u)f(u)\Big)\p(v)\ol f(v)
\\-\sum_{v\in V_*}
\Big(\p(v)-\sum_{\be=(v,\,u)\in\cA_*}
\frac{\p(u)}{\sqrt{\vk_v\vk_u}}\Big)\p(v)|f(v)|^2
\\=\sum_{\be=(v,\,u)\in\cA_*}
c_{uv}\big(|f(v)|^2-
e^{i\lan\t(\be),\,\vt\rangle}f(u)\ol f(v)\big),
\end{multline}
where $c_{uv}$ is defined in \eqref{eq.L1}. By the identity $\t(\be)=-\t(\ul\be)$
for all $\be\in\cA_*$, we have
\begin{multline}
\lb{eq.L23}
\sum_{\be=(v,\,u)\in\cA_*}
c_{uv}\big(|f(v)|^2-
e^{i\lan\t(\be),\,\vt\rangle}f(u)\ol f(v)\big)\\=\sum_{\ul\be=(u,\,v)\in\cA_*}
c_{uv}\big(|f(v)|^2-
e^{-i\lan\t(\ul\be),\,\vt\rangle}f(u)\ol f(v)\big)\\
=\sum_{\be=(v,\,u)\in\cA_*}
c_{uv}\big(|f(u)|^2-
e^{-i\lan\t(\be),\,\vt\rangle}f(v)\ol f(u)\big).
\end{multline}
Then the identity \eqref{eq.L22} can be rewritten in the form:
\begin{multline}
\lb{eq.L24}
\big\lan\big(H(\vt)-\l_1(0)\1\big)\p f,\p f\big\rangle_{V_*}
\\
={\frac12}\sum_{\be=(v,\,u)\in\cA_*}
c_{uv}\big(|f(v)|^2-
e^{i\lan\t(\be),\,\vt\rangle}f(u)\ol f(v)
+|f(u)|^2-e^{-i\lan\t(\be),\,\vt\rangle}f(v)\ol f(u)\big)
\\
={\frac12}\sum_{\be=(v,\,u)\in\cA_*}
c_{uv}\big|f(v)-e^{i\lan\t (\be),\,\vt\rangle}f(u)\big|^2,
\end{multline}
as required.

\emph{ii}) Let $g=g(\vt)\in\ell^2(V_*)$ be an eigenvector of the operator $$H(\vt)-\l_1(0)\1,$$
corresponding to the smallest eigenvalue $\l_1(\vt)-\l_1(0)$, such that
$\|\vk^{1/2}f_0\|_{V_*}=1$, where $f_0=f_0(\vt)=\p^{-1}g$, and $\vk$ is the operator
of multiplication by the function $\vk(v)=\vk_v$. Here $\|\,\cdot\,\|_{V_*}$
is the norm in the space $\ell^2(V_*)$. Then, by the minimax principle (see
the matrix property i) and the identities \eqref{qflo}, \eqref{eq.L1}, we have
$$
\begin{aligned}
\l_1^{(0)}(\vt)&\!=\!\!\min\limits_{f\in\ell^2(V_*) \atop \|f\|_{V_*}=1}\langle\D(\vt)f,f\rangle_{V_*}
\!\leq\!
\langle\D(\vt)\vk^{1/2}f_0,\vk^{1/2}f_0\rangle_{V_*}
\\
&\!=\!
\frac12\!\!\sum_{\be=(v,\,u)\in\cA_*}\!\!\big|f_0(v)\!-
e^{i\lan\t(\be),\vt\rangle}\,f_0(u)\big|^2
\\
&\!\leq\!\frac1{2\p_-^2}\!\!\sum_{\be=(v,\,u)\in\cA_*}\!\!
c_{uv}\big|f_0(v)-
e^{i\lan\t(\be),\vt\rangle}\,f_0(u)\big|^2
\\
&\!=\!
{1\over \p_-^2}\lan\big(H(\vt)-\l_1(0)\1\big)\p f_0,\p f_0\rangle_{V_*}
\\&\!=\!
{\|\p f_0\|_{V_*}^2\over \p_-^2}\,
\big(\l_1(\vt)\!-\!\l_1(0)\big)
\!\leq\! c_0^2\big(\l_1(\vt)\!-\!\l_1(0)\big), \quad \forall\,\vt\!\in\!\T^d.
\end{aligned}
$$
Thus, the first inequality in \eqref{eq.T1} is proved.

Similarly, let $g_0=g_0(\vt)\in\ell^2(V_*)$ be an eigenvector of the operator
$\D(\vt)$, corresponding to the smallest eigenvalue $\l_1^{(0)}(\vt)$, such that
$\|\vk^{-1/2}f_0\|_{V_*}=1$, where $f_0=f_0(\vt)=\p g_0$. Then
\begin{multline*}
\l_1(\vt)-\l_1(0)=\min\limits_{f\in\ell^2(V_*) \atop \|f\|_{V_*}=1}\langle\big(H(\vt)-\l_1(0)\1\big)f,f\rangle_{V_*}
\\
\leq
\langle\big(H(\vt)-\l_1(0)\1\big)\vk^{-1/2}f_0,\vk^{-1/2}f_0\rangle_{V_*}
\\
=\langle\big(H(\vt)-\l_1(0)\1\big)\p\p^{-1}\vk^{-1/2}f_0,
\p\p^{-1}\vk^{-1/2}f_0\rangle_{V_*}
\\
={1\over 2}\sum_{\be=(v,\,u)\in\cA_*}
c_{uv}\bigg|{f_0(v)\over \p(v)\sqrt{\vk_v}}-e^{i\lan\t (\be),\,\vt\rangle}{f_0(u)\over \p(u)\sqrt{\vk_u}}\bigg|^2
\\
\leq{\psi_+^2\over 2}\sum_{\be=(v,\,u)\in\cA_*}
\bigg|{f_0(v)\over \p(v)\sqrt{\vk_v}}-e^{i\lan\t (\be),\,\vt\rangle}{f_0(u)\over \p(u)\sqrt{\vk_u}}\bigg|^2
\\
=\!\psi_+^2\lan\D(\vt)\p^{-1}f_0,\p^{-1}f_0\rangle_{V_*}
\!=\!\psi_+^2\l_1^{(0)}(\vt)\|\p^{-1}f_0\|_{V_*}^2\!\leq\!
c_0^2\,\l_1^{(0)}(\vt), \quad \forall\,\vt\in\T^d,
\end{multline*}
and the second inequality in \eqref{eq.T1} is proved. \qq $\BBox$

\medskip

\no \textbf{Proof of Theorem \ref{TMR1}.} \emph{ii}) Using \eqref{ban.1}, \eqref{bes} and
\eqref{eq.T1}, for some $\vt_+\in\T^d$ we have
\begin{equation}
\lb{qqq1}
|\s_1(H)|=\l_1^+-\l_1(0)=\l_1(\vt_+)-\l_1(0)\leq
c_0^{2}\,\l_1^{(0)}(\vt_+)
\leq c_0^{2}\,\l_1^{0+}=c_0^{2}\,|\s_1(\D)|.
\end{equation}
Similarly, for some $\vt_+^0\in\T^d$ we have
\begin{equation}
\lb{qqq2}
\begin{aligned}
c_0^{-2}|\s_1(\D)|=
c_0^{-2}\l_1^{0+}&=
c_0^{-2}\l_1^{(0)}(\vt_+^0)
\\
&\leq\l_1(\vt_+^0)-\l_1(0)
\leq\l_1^+-\l_1(0)=|\s_1(H)|.
\end{aligned}
\end{equation}
The formulas \eqref{eq.TMR1} follow from \eqref{qqq1} and \eqref{qqq2}.

We recall the known fact that the first spectral band
$
\s_1(\D)=[0,\l_1^{0+}]
$
of the Laplacian~$\D$ is non-degenerate. Indeed, assume that $0$ is an eigenvalue
of the operator $\D$ with infinite multiplicity. Then (see, e.g., Theorem 4.5.2 in
\cite{BK13}) there exists an eigenfunction $0\neq f\in\ell^2(V)$, corresponding to
this eigenvalue and having finite support. By \eqref{qflo1}, we obtain
$$
0=\lan \D f,f\rangle_V={1\over 2}\sum_{(v,\,u)\in\cA}
\bigg|{f(v)\over \sqrt{\vk_v}}-{f(u)\over \sqrt{\vk_u}}\,\bigg|^2,
$$
which yields
$$
{f(u)\over \sqrt{\vk_u}}={f(v)\over \sqrt{\vk_v}}, \quad \forall\,u,v\in V
\textrm{ such that } u\sim v.
$$
Since the support of the function $f$ is finite and the graph $\G$ is connected,
we conclude that $f=0$. We get a contradiction. Thus, the first spectral band
$\s_1(\D)$ of the
Laplacian~$\D$ is non-degenerate. Then, using the fact that $\p$ is a vector with
positive components, from the inequality \eqref{eq.TMR1} we conclude that the
first spectral band $\s_1(H)$ of the Schr\"odinger operator $H=\D+Q$ is non-degenerate.

\emph{iii}) The effective mass tensor $m=M^{-1}$, where $M=\{M_{ij}\}$ is defined in
\eqref{tem}. Then the estimate \eqref{eem} follows from the inequality \eqref{eq.T1}. \qq $\BBox$

\subsection{Estimate of the Lebesgue measure of the spectrum}
We need the following representation of the fiber Schr\"odinger operator $H(\vt)$,
$\vt\in\T^d$:
\begin{equation}
\label{eq.1} H(\vt)=H_0+h(\vt),\quad
H_0={1\over (2\pi)^d}\int\limits_{\T^d}H(\vt )d\vt.
\end{equation}
Using the formulas \eqref{eq.1} and \eqref{Devt}, we obtain the representation
of the operator
$$
h(\vt)=\{h_{uv}(\vt)\}_{u,v\in V_*}
$$
in the standard basis \eqref{sonb}:
\begin{equation}
\label{tl2.15}
h_{uv}(\vt)={-1\over \sqrt{\vk_u\vk_v}}\sum\limits_{\be=(u,v)\in\cB_*}e^{-i\lan\t
(\be),\,\vt\rangle},
\end{equation}
where $\cB_*$ is the set of the fundamental graph edges with non-zero indices.

\medskip

\no \textbf{Proof of Theorem \ref{T2.2}.}
\emph{i}) We define the diagonal matrix $B(\vt)$ in the following way:
\begin{equation}
\lb{eq.2'}
B(\vt)=\diag(B_u(\vt))_{u\in V_*},\quad B_u(\vt)=\sum_{v\in V_*}|h_{uv}(\vt)|, \quad \vt\in\T^d.
\end{equation}
From \eqref{tl2.15} it follows that
\begin{equation}
\lb{wvv1} |h_{uv}(\vt)|\leq|h_{uv}(0)|={\gz_{uv}\over \sqrt{\vk_u\vk_v}}, \quad \forall\,(u,v,\vt)\in V_*^2\ts\T^d,
\end{equation}
where $\gz_{uv}$ is the number of edges of the form $(u,v)$ with non-zero indices on the
fundamental graph $\G_*$. From \eqref{eq.2'}, \eqref{wvv1} we deduce that
\begin{equation}
\lb{B0}
B(\vt)\leq B(0),\quad\forall\,\vt\in\T^d.
\end{equation}
The estimate \eqref{B0} and the matrix property vii) give
\begin{equation}
\lb{eq.2} -B(0) \leq -B(\vt)\le h(\vt)\le B(\vt)\leq B(0),\quad
\forall\,\vt\in\T^d.
\end{equation}
Combining \eqref{eq.1} with \eqref{eq.2}, we obtain
$$
H_0-B(0)\le H(\vt)\le H_0+B(0),
$$
whence
$$
\l_n(H_0-B(0))\leq\l_n^-\le \l_n(\vt)\leq\l_n^+\le \l_n(H_0+B(0)),
\qqq \forall\,(n,\vt)\in\N_\n\ts\T^d,
$$
which yields
\begin{equation}
\lb{009}
\big|\s(H)\big|\le \sum_{n=1}^{\nu}(\l_n^+-\l_n^-)
\leq\sum_{n=1}^\nu\big(\l_n(H_0+B(0))-\l_n(H_0-B(0))\big)=
 2\Tr B(0).
\end{equation}
Using the relations \eqref{eq.2'}, \eqref{wvv1}, we have
\begin{equation}
\lb{pro} 2\Tr B(0)=2\sum_{u\in V_*}B_u(0)=2\sum_{u,v\in V_*}|h_{uv}(0)|=
2\sum_{u,v\in V_*}{\gz_{uv}\over \sqrt{\vk_u\vk_v}}\,.
\end{equation}
By the Cauchy-Schwarz inequality, we obtain
\begin{equation}
\lb{CS}
\sum_{u,v\in V_*}{\gz_{uv}\over \sqrt{\vk_u\vk_v}}\leq
\bigg(\sum_{u,v\in V_*}{\gz_{uv}\over \vk_u}\bigg)^{1/2}
\bigg(\sum_{u,v\in V_*}{\gz_{uv}\over \vk_v}\bigg)^{1/2}=
\sum_{v\in V_*}{\gz_v\over \vk_v}\,.
\end{equation}
Here we have used the identities $\gz_{uv}=\gz_{vu}$ and $\gz_v=\sum_{u\in V_*}\gz_{vu}$
for all $u,v\in V_*$. The estimate \eqref{eq.7} follows from \eqref{009}--\eqref{CS}.

Now, we prove \eqref{GEga}. Since $\l_1^-$ and $\l_\n^+$ are the lower and upper
endpoints of the spectrum $\s(H)$, using the estimate \eqref{eq.7}, we obtain
\begin{equation}
\lb{GEga1}
\sum_{n=1}^s|\g_n(H)|=
\l^+_\n-\l_1^--\big|\s(H)\big|\geq\l^+_\n-\l_1^--2\gz.
\end{equation}
We rewrite the sequence $(Q(v))_{v\in V_*}$ in nondecreasing order
\begin{equation}
\lb{wtqn}
q_1^\bu\le q_2^\bu \le\ldots \le q_\n^\bu, \quad  q_1^\bu=0,
\end{equation}
where $q_1^\bu=Q(v_1)$, $q_2^\bu=Q(v_2),\ldots,q_\n^\bu=Q(v_\n)$
for some distinct vertices $v_1,v_2,\ldots,v_\n\in V_*$, and without
loss of generality we may assume that $q_1^\bu=0$.

Then, according to the matrix property ii), the eigenvalues of the fiber operator
$H(\vt)=\D(\vt)+Q$ satisfy the inequalities
\begin{equation}
\lb{qq0}
\begin{aligned}
&q_n^\bu\le q_n^\bu+\l_1^0(\vt)\le\l_n(\vt)\le q_n^\bu+ \l_\n^0(\vt)\le q_n^\bu+\l_\n^{0+},
  \\
&\l_n^0(\vt)\le\l_n(\vt)\le \l_n^0(\vt)+q_\n^\bu, \quad \forall\, (\vt,n)\in\T^d\ts\N_\n.
\end{aligned}
\end{equation}
From the first inequality in \eqref{qq0} we obtain
\begin{equation}
\lb{qq}
\l_\n^+\geq q_\n^\bu,\quad \l_1^-\leq \l_\n^{0+},
\end{equation}
and, using the second inequality in \eqref{qq0}, we have
\begin{equation}
\lb{qq1}
\l_\n^{0+}=\max_{\vt\in\T^d}\l_\n^{0}(\vt)=
\l_\n^{0}(\vt_+)\leq\l_\n(\vt_+)\leq \l_\n^+,
\end{equation}
\begin{equation}
\lb{qq2}
0=\l_1^{0-}=\min_{\vt\in\T^d}\l_1^{0}(\vt)=\l_1^{0}(\vt_-)\geq
\l_1(\vt_-)-q_\n^\bu\geq\l_1^--q_\n^\bu
\end{equation}
for some $\vt_-,\vt_+\in\T^d$. From \eqref{qq}--\eqref{qq2} it follows that
$$
\l^+_\n-\l_1^-\geq q_\n^\bu-\l_\n^{0+},\quad
\l^+_\n-\l_1^-\geq \l_\n^{0+}-q_\n^\bu,
$$
which yields $\l^+_\n-\l_1^-\geq C_0$, where $C_0$ is defined in \eqref{GEga}.
Thus, the estimate~\eqref{GEga} is proved.

The proof of item ii) will be given in the proof of Theorem~\ref{T100}.ii.
\qq $\BBox$

\section{\lb{Sec5} Loop and bipartite graphs}
\setcounter{equation}{0}
\subsection{The Schr\"odinger operator on loop graphs}

\no \textbf{Proof of Theorem \ref{T100}.} \emph{i}) Let $\G$ be a loop graph. Then, by
\eqref{eq.1}, \eqref{tl2.15}, we obtain
$$
H(\vt)=H_0+h(\vt),
$$
where the operator $h(\vt)$, $\vt\in\T^d$, in the standard basis \eqref{sonb}
has the form
\begin{equation}
\lb{ele}
\begin{aligned}
&h(\vt)=\diag\big(h_{vv}(\vt)\big)_{v\in V_*},
\\
&h_{vv}(\vt)=-{1\over \vk_v}\sum\limits_{\be=(v,v)\in\cB_*}
\cos\lan\t(\be),\vt\rangle,
\end{aligned}
\end{equation}
which yields $h_{vv}(0)\leq h_{vv}(\vt)$ for all $(\vt,v)\in\T^d\ts V_*$. Then
$h(0)\leq h(\vt)$ and we have
$$
H(0)=H_0+h(0)\leq H_0+h(\vt)=H(\vt), \quad \forall\,\vt\in\T^d,
$$
whence
$\l_n(0)\leq\l_n(\vt)$ for all $(\vt,n)\in\T^d\ts\N_\n$, which yields
$$
\l_n^-=\min\limits_{\vt\in\T^d}\l_n(\vt)=\l_n(0).
$$

\emph{ii}) Let $\G$ be a precise loop graph with a precise quasimomentum $\vt_0\in\T^d$. Then
from the identity \eqref{ele} it follows that $h_{vv}(\vt)\leq h_{vv}(\vt_0)$
for all $(\vt,v)\in\T^d\ts V_*$, since $\cos\lan\t(\be),\vt_0\rangle=-1$
for all edges $\be\in\cB_*$ with non-zero indices. Consequently, $h(\vt)\le h(\vt_0)$
and
$$
H(\vt)=H_0+h(\vt)\leq H_0+h(\vt_0)=H(\vt_0),
$$
whence $\l_n(\vt)\le\l_n(\vt_0)$ for all $(\vt,n)\in\T^d\ts\N_\n$.
This yields
$$
\l_n^+=\max\limits_{\vt\in\T^d}\l_n(\vt)=\l_n(\vt_0)
$$
and, by \eqref{eq.5}, the spectral bands $\s_n$ have the form~\eqref{es}.

Using the formulas \eqref{es} and \eqref{ele}, we obtain
$$
\sum\limits_{n=1}^\n|\s_n|=\sum\limits_{n=1}^\n\big(\l_n(\vt_0)-\l_n(0)\big)
=\Tr\big(H(\vt_0)-H(0)\big)=\Tr\big(h(\vt_0)-h(0)\big)=
2\sum_{v\in V_*}\frac{\gz_v}{\vk_v},
$$
and the identity \eqref{es.1} is proved.

Now, let all edges of the fundamental graph $\G_*$ with non-zero indices have the
form $(v,v)$ for some vertex $v\in V_*$. It only remains to prove the first
identity in \eqref{loo}. Without loss of generality we may assume that the operator
$H(\vt)$ in the standard basis \eqref{sonb} has the form
$$
H(\vt )=\left(
\begin{array}{cc}
  A & y \\
  y^* & \D_{vv}(\vt)+Q(v)
\end{array}\right),
$$
where $y\in\C^{\nu-1}$, the entry $\D_{vv}(\vt)$ is given by the formula \eqref{Devt},
and $A$ is a self-adjoint $(\nu-1)\ts(\nu-1)$-matrix not depending on $\vt$. The
eigenvalues $\m_1\leq\ldots\leq\m_{\n-1}$ of the matrix $A$ do not depend on $\vt$.
Then, applying the matrix property iii), we obtain
$$
\l_1(\vt)\leq\mu_1\le\l_2(\vt)\leq\m_2\leq\dots\leq
\mu_{\n-1}\le\l_{\n}(\vt),\quad \forall\,\vt\in\T^d.
$$
This yields that the spectral bands of the operator $H$ may only touch, but do not
overlap, i.e.,
$$
|\s(H)|=\sum\limits_{n=1}^\n|\s_n|.
$$
Thus, in this case the estimate \eqref{eq.7} and, consequently, the first inequality
in \eqref{GEga} become identities. \qq $\BBox$

\medskip

\no \textbf{Remark.} If there exists $\vt_0\in\T^d$ such that $\lan\t(\be),\vt_0 \rangle/\pi$
is odd for all edges $\be\in\cB_*$ with non-zero indices, then $\vt_0$ is a precise
quasimomentum for $\G$.

\subsection{Laplacians on bipartite graphs} A graph is called \emph{bipartite}
if its vertex set is divided into two disjoint sets (called \emph{parts} of the graph)
such that each edge of the graph connects vertices from the distinct parts.
Examples of bipartite graphs are the $d$-dimensional lattice and the hexagonal lattice.
The face-centered cubic lattice (Fig.~\ref{ff.FCC}) is non-bipartite. It is known
(see, e.g., \cite{MW89}), that the following statements are equivalent:

$\bu$ a graph $\G$ is bipartite;

$\bu$ the point $2$ belongs to the spectrum $\s(\D)$ of the Laplacian $\D$ on the
graph~$\G$;

$\bu$ the spectrum $\s(\D)$ is symmetric with respect to the point 1.

We formulate some spectral properties of the Laplacian on bipartite periodic graphs.

\begin{theorem} \lb{TBG}
The following statements hold.

i) The fundamental graph $\G_*$ is bipartite if and only if for each $\vt\in\T^d$
the spectrum of the fiber operator $\D(\vt)$ is symmetric with respect to the point~$1$.

ii) Let the fundamental graph $\G_*$ be bipartite with parts $V_1,V_2$ and
$m=\# V_1-\# V_2>0$. Then $\{1\}$ is a flat band of the Laplacian $\D$ on
the periodic graph $\G$ of multiplicity at least $m$.

iii) If there exist $s$ gaps $\g_1(\D),\ldots,\g_s(\D)$ in the spectrum of
the Laplacian $\D$ on a bipartite periodic graph $\G$, then the following estimate
holds true:
\begin{equation}
\lb{GEgaB}
\sum_{n=1}^s|\g_n(\D)|\ge
2(1-\gz),
\end{equation}
where $\gz$ is defined in \eqref{eq.7}.

iv) Let $\G$ be a bipartite loop graph (the fundamental graph $\G_*$
is non-bipartite, since all edges of the graph $\G_*$ with non-zero indices are loops). Then each spectral band of the Laplacian $\D$ on $\G$ has the form
$$
\s_n^0=\big[\l_n^{0-},\l_n^{0+}\big],\quad n\in \N_\n,
$$
where $\l_n^{0-}$ and $\l_n^{0+}$ are the eigenvalues of the operators $\D(0)$ and
$2\1-\D(0)$, respectively ($\1$ is the identity operator).
\end{theorem}

\no \textbf{Proof.} \emph{ii}) Let $\n_s=\# V_s$, $s=1,2$. Since vertices from the same part of
the bipartite graph $\G_*=(V_*,\cE_*)$ are not adjacent to each other, for each
$\vt \in\T^d$ the fiber Laplacian \eqref{l2.13} in the standard basis \eqref{sonb}
can be represented in the following form:
\begin{equation}
\label{m1}
\D(\vt)=\begin{pmatrix}
  \1_{\n_1} & A(\vt)\\[10pt]
  A^*(\vt) & \1_{\nu_2}
\end{pmatrix},
\end{equation}
where $\1_n$ is the identity $(n\ts n)$-matrix, and $A(\vt)$ is some $(\n_1\ts\n_2)$-matrix.
Then $\l=1$ is an eigenvalue of the matrix $\D(\vt)$ with an eigenfunction
$$
f=(f_1,f_2)\in\ell^2(V_*),
$$
 where $f_s\in\ell^2(V_s)$, if and only if
\begin{equation}
\lb{fst}
A(\vt)f_2=0,\quad A^*(\vt)f_1=0.
\end{equation}
The system $A(\vt)f_2=0$ always has a zero solution $f_2=0$. Since
$\rank A^*(\vt)\leq\n_2$, the number of linear independent solutions of the system
$$
A^*(\vt)f_1=0
$$
 is equal to
$$
\n_1-\rank A^*(\vt)\geq\n_1-\n_2.
$$
Thus, for each  $\vt\in\T^d$ the operator $\D(\vt)$ has the eigenvalue $\l=1$ of
multiplicity at least $\n_1-\n_2$, which yields the required statement.

The proofs of the remaining items repeat the arguments for the combinatorial
Laplacian (see~\cite[the proof of Theorem~5.2]{KS14}). \qq $\BBox$

\section{\lb{Sec6} Perturbations of the $d$-dimensional lattice}
\setcounter{equation}{0}

We consider $d$-dimensional lattice $\dL^d=(V,\cE)$, where the vertex set and
the edge set are given by
\begin{equation}
\lb{dLg} V=\Z^d,\qquad
\cE=\big\{(\mm,\mm+a_s), \quad
\forall\,\mm\in\Z^d, \qq s\in\N_d\big\},
\end{equation}
and the orthonormal basis $a_1,\ldots,a_d$ coincides with the periods of the lattice
$\dL^d$. The minimal fundamental graph $\dL^d_*$ of the lattice $\dL^d$ consists
of one vertex $v$ and $2d$ oriented loop edges
$$
\be_1=\ul\be_1=\ldots=\be_d=\ul\be_d=(v,v)
$$
with indices $\pm a_1,\ldots,\pm a_d$. For $\vt_\pi=(\pi,\ldots,\pi)\in\T^d$ we have
the identity
$$
\cos\lan\t(\be),\vt_\pi\rangle=-1
$$
for all edges of the fundamental graph $\dL^d_*$. Consequently, the graph $\dL^d$
is a precise loop  graph with the precise quasimomentum $\vt_\pi$. It is known that
the spectrum of the Laplacian $\D$ on $\dL^d$ has the form
$$
\s(\D)=\s_{ac}(\D)=[0,2].
$$

\setlength{\unitlength}{1.0mm}
\begin{figure}[h]
\centering
\unitlength 1mm 
\linethickness{0.4pt}
\ifx\plotpoint\undefined\newsavebox{\plotpoint}\fi 
\begin{picture}(120,50)(0,0)

\put(10,10){\line(1,0){40.00}} \put(10,30){\line(1,0){40.00}}
\put(10,50){\line(1,0){40.00}} \put(10,10){\line(0,1){40.00}}
\put(30,10){\line(0,1){40.00}} \put(50,10){\line(0,1){40.00}}

\put(10,10){\circle{1}} \put(30,10){\circle{1}}
\put(50,10){\circle{1}}

\put(10,30){\circle{1}} \put(30,30){\circle{1}}
\put(50,30){\circle{1}}

\put(10,50){\circle{1}} \put(30,50){\circle{1}}
\put(50,50){\circle{1}}

\put(10,10){\vector(1,0){20.00}} \put(10,10){\vector(0,1){20.00}}

\put(7.5,7.5){$\scriptstyle v_\n$} \put(1,29){$\scriptstyle
v_\n+a_2$} \put(29,7){$\scriptstyle v_\n+a_1$}
\put(20,8){$\scriptstyle a_1$} \put(6,20){$\scriptstyle a_2$}
\put(14,26.5){$\scriptstyle v_1$} \put(23.8,17){$\scriptstyle
v_{\nu-2}$}

\put(24.0,12){$\scriptstyle v_{\nu-1}$} \put(20,26.5){$\scriptstyle
v_2$}

\put(10,10){\line(2,3){10.00}} \put(10,10){\line(3,2){15.00}}
\put(10,10){\line(1,3){5.00}} \put(10,10){\line(3,1){15.00}}
\put(20.5,21){$\ddots$} \put(15,25){\circle{1}}
\put(20,25){\circle{1}} \put(25,15){\circle{1}}
\put(25,20){\circle{1}}

\put(30,10){\line(2,3){10.00}} \put(30,10){\line(3,2){15.00}}
\put(30,10){\line(1,3){5.00}} \put(30,10){\line(3,1){15.00}}
\put(40.5,21){$\ddots$} \put(35,25){\circle{1}}
\put(40,25){\circle{1}} \put(45,15){\circle{1}}
\put(45,20){\circle{1}}

\put(10,30){\line(2,3){10.00}} \put(10,30){\line(3,2){15.00}}
\put(10,30){\line(1,3){5.00}} \put(10,30){\line(3,1){15.00}}
\put(20.5,41){$\ddots$} \put(15,45){\circle{1}}
\put(20,45){\circle{1}} \put(25,35){\circle{1}}
\put(25,40){\circle{1}}

\put(30,30){\line(2,3){10.00}} \put(30,30){\line(3,2){15.00}}
\put(30,30){\line(1,3){5.00}} \put(30,30){\line(3,1){15.00}}
\put(40.5,41){$\ddots$} \put(35,45){\circle{1}}
\put(40,45){\circle{1}} \put(45,35){\circle{1}}
\put(45,40){\circle{1}} \put(-4,10){({a})}
\put(64,25){({b})}

\multiput(80,45)(4,0){5}{\line(1,0){2}}
\multiput(100,25)(0,4){5}{\line(0,1){2}}
\put(80,25){\line(1,0){20.00}} \put(80,25){\line(0,1){20.00}}
\put(80,25){\circle{1}} \put(100,25){\circle{1}}

\put(80,45){\circle{1}} \put(100,45){\circle{1}}

\put(80,25){\vector(1,0){20.00}} \put(80,25){\vector(0,1){20.00}}

\put(77.0,23.0){$\scriptstyle v_\n$} \put(76.5,45){$\scriptstyle
v_\n$} \put(101,23){$\scriptstyle v_\n$} \put(101,45){$\scriptstyle
v_\n$} \put(90,23){$\scriptstyle a_1$} \put(76,35){$\scriptstyle
a_2$}

\put(84,41.5){$\scriptstyle v_1$} \put(93.8,32.5){$\scriptstyle
v_{\nu-2}$}

\put(94.0,27){$\scriptstyle v_{\nu-1}$} \put(90,41.5){$\scriptstyle
v_2$}

\put(80,25){\line(2,3){10.00}} \put(80,25){\line(3,2){15.00}}
\put(80,25){\line(1,3){5.00}} \put(80,25){\line(3,1){15.00}}
\put(90.5,36){$\ddots$} \put(85,40){\circle{1}}
\put(90,40){\circle{1}} \put(95,30){\circle{1}}
\put(95,35){\circle{1}}
\put(64,10){({c})}

\put(80,10){\line(1,0){30.00}} \put(80,9){\line(0,1){2.00}}
\put(110,9){\line(0,1){2.00}} \put(90,9){\line(0,1){2.00}}
\put(100,9){\line(0,1){2.00}} \put(95,10){\circle*{1}}

\put(80,9.8){\line(1,0){10.00}} \put(80,10.2){\line(1,0){10.00}}

\put(100,9.8){\line(1,0){10.00}} \put(100,10.2){\line(1,0){10.00}}

\put(103,12){$\s_\n^0$} \put(83,12){$\s_1^0$} 
\put(79,6){$\scriptstyle0$} \put(94.3,6){$\scriptstyle1$}
\put(109.3,6){$\scriptstyle2$}
\end{picture}
\vspace{-0.5cm} \caption{\footnotesize {a}) A $\Z^d$-periodic
graph $\G$;\quad {b}) the fundamental graph $\G_*$; \quad
{c})~the spectrum of the Laplacian.} \label{ff.10}
\end{figure}

We describe the spectrum of the Laplace and Schr\"odinger operators on perturbations
$\G$ of the lattice $\dL^d$, shown in Fig. \ref{ff.10}\emph{a}; these perturbations are precise
loop graphs.

\begin{proposition}\lb{TG1}
Let $\G_*$ be obtained from the fundamental graph $\dL^d_*$ of
$d$\nobreakdash-di\-men\-sional lattice $\dL^d$ by adding $\n-1\geq1$ vertices
$v_1,\ldots,v_{\n-1}$ and $\n-1$ unoriented edges
$(v_1,v_\n),\ldots,(v_{\n-1},v_\n)$ with zero indices {\rm (}see
Fig. {\rm \ref{ff.10}b)}, $v_\n$ is a unique vertex of $\dL^d_*$. Then $\G$
is a precise loop graph with the precise quasimomentum $\vt_\pi=(\pi,\ldots,\pi)\in\T^d$
and the following statements hold.

i) The spectrum of the Laplacian $\D$ on the graph $\G$ has the form
\begin{equation}
\lb{acs1} \s(\D)=\s_{ac}(\D)\cup\s_{fb}(\D),\quad \s_{fb}(\D)=\{1\},
\end{equation}
where the flat band $\{1\}$ has multiplicity $\n-2$, and the absolutely continuous spectrum
$\s_{ac}(\D)$ consists of two bands $\s_1^0$ and $\s_\n^0:$
\begin{equation}
\lb{acs2} \s_{ac}(\D)=\s_1^0\cup\s_\n^0, \quad
\s_1^0=[0,{\textstyle\frac{2d}{\x}}],\quad
\s_\n^0=[{\textstyle2-\frac{2d}{\x}},2],\quad \x=\n-1+2d.
\end{equation}

ii) The spectrum of Schr\"odinger operator $H=\D+Q$ on the graph $\G$
has the form
\begin{equation}
\lb{acs333} \s(H)=\bigcup_{n=1}^\n\big[\l_n(0),\l_n(\vt_\pi)\big],\quad
\vt_\pi=(\pi,\ldots,\pi).
\end{equation}

iii) Let $q_j=Q(v_j)$, $j\in\N_\n$. Without loss of generality we may assume
that $q_\n=0$. If all values $q_1,\ldots,q_{\n-1}$ of the potential at the vertices
of the fundamental graph $\G_*$ are distinct, then $\s(H)=\s_{ac}(H)$, i.e.,
$\s_{fb}(H)=\varnothing $.

iv) Let among the numbers $q_1,\ldots,q_{\n-1}$ there exist a value $q_*$
of multiplicity~$m$. Then the spectrum of the Schr\"odinger operator $H$ on $\G$
has the flat band $\{q_*+1\}$ of multiplicity $m-1$.

v) The Lebesgue measure of the spectrum of the Schr\"odinger operator $H$ on
$\G$ satisfies
\begin{equation}
\lb{sp2} \textstyle|\s(H)|=\frac{4d}{\x}.
\end{equation}
\end{proposition}
\no \textbf{Proof.} \emph{i}) -- \emph{ii}) The fundamental graph $\G_*$ consists of $\n\geq2$ vertices
$v_1,\ldots,v_\n$; ${\n-1}$ edges ${(v_1,v_\n),\ldots,(v_{\n-1},v_\n)}$ with zero
indices and $2d$ oriented loops at the vertex $v_\n$ with the indices
$\pm a_1,\ldots,\pm a_d$. Since all edges of the fundamental graph $\G_*$
with non-zero indices are loops and $\cos\lan\t(\be),\vt_\pi\rangle=-1$ for all
such edges $\be\in\cB_*$, the graph $\G$ is a precise loop graph with
the precise quasimomentum $\vt_\pi=(\pi,\ldots,\pi)$. Then, by Theorem \ref{T100}.ii,
the spectral bands of the Schr\"odinger operator $H$ are given by
$$
\s_n=[\l_n^-,\l_n^+]=[\l_n(0),\l_n(\vt_\pi)],\quad\forall\, n\in\N_\n,
$$
and item \emph{ii}) is proved.

By (\ref{Devt}), for each $\vt=(\vt_j)_{j=1}^d\in\T^d$ we have
\begin{equation}
\lb{fc0}
\D(\vt)=\left(
\begin{array}{cccc}
  1 & 0  & \ldots & -\frac1{\sqrt{\x}}  \\[2pt]
  0 & 1 & \ldots & -\frac1{\sqrt{\x}} \\[2pt]
  \ldots & \ldots & \ldots & \ldots  \\[2pt]
  -\frac1{\sqrt{\x}} & -\frac1{\sqrt{\x}} & \ldots & 1-\frac{2c_0}{\x}  \\
\end{array}\right),\quad
\ca c_0=c_1+\ldots+c_d\\
c_j=\cos\vt_j, \qq j\in\N_d\\
\x=\n-1+2d\ac\!\!.
\end{equation}
Using the formula \eqref{det}, we obtain
$$
\begin{aligned}
\det\big(\D(0)-\l \1_\n\big)&=(1-\l)^{\n-2}\,\l\,
\Big(\l-2+{2d\over \x}\Big),
\\
\det\big(\D(\vt_\pi)-\l \1_\n\big)&=(1-\l)^{\n-2}\,(\l-2)\,
\Big(\l-{2d\over \x}\Big),
\end{aligned}
$$
where $\1_\n$ is the identity $(\n\ts\n)$-matrix. Then the eigenvalues of the matrices
$\D(0)$ and $\D(\vt_\pi)$ have the form
$$
\textstyle\l_1^0(0)=0,\quad \l_2^0(0)=\ldots=\l_{\n-1}^0(0)=1,\quad
\l_\n^0(0)=2-{2d\over \x},
$$
$$
\textstyle\l_1^0(\vt_\pi)={2d\over \x},\quad \l_2^0(\vt_\pi)=\ldots=
\l_{\n-1}^0(\vt_\pi)=1,\quad \l_\n^0(\vt_\pi)=2.
$$
Thus, by \eqref{acs333}, the spectrum of the Laplacian on the graph $\G$
satisfies the identities \eqref{acs1}, \eqref{acs2}.

\emph{iii}) For each $\vt\in\T^d$ the operator $H(\vt)$ has the form
$$
H(\vt)=\D(\vt)+q,\quad q=\diag\,(q_1,\ldots,q_{\n-1},0),
$$
where $\D(\vt)$ is given by \eqref{fc0}. Using the formula \eqref{det},
we write the characteristic polynomial of the matrix $H(\vt)$ in the following form:
\begin{equation}
\lb{lco1}
\det(H(\vt)-\l\1_\n)=\textstyle\frac1{\x}\,
\big[(\x-\x\l-2c_0)W(\l)+W'(\l) \big],
\end{equation}
where
\begin{equation}
\lb{WWW}
W(\l)=(1+q_1-\l)\ldots(1+q_{\n-1}-\l).
\end{equation}

We show that $\s(H)=\s_{ac}(H)$ by contradiction. Suppose that the Schr\"odinger
operator $H$ has an eigenvalue $\l$ of infinite multiplicity. Then
$\det(H(\vt)-\l\1_\n)=0$ for all $\vt\in\T^d$. The linear combination \eqref{lco1}
of the linearly independent functions $c_1,\ldots,c_d,1$ is identically equal to 0
if and only if
\begin{equation}
\lb{sys} W(\l)=0,\quad  W'(\l)=0.
\end{equation}
All values $q_1,\ldots,q_{\n-1}$ of the potential are distinct. Therefore,
each zero of the function $W$ is simple. This contradicts the identities \eqref{sys}.
Consequently, $\s(H)=\s_{ac}(H)$.

\emph{iv})  Without loss of generality we may assume that
$$
q_1=q_2=\ldots=q_m=q_*.
$$
Then $\l=q_*+1$ is a zero of multiplicity $m$ of the function $W$, defined by the
formula \eqref{WWW}, and is a zero of multiplicity $m-1$ of the function $W'$.
Consequently, the system \eqref{sys}, which defines all eigenvalues of infinite
multiplicity, has the solution $\l=q_*+1$ of multiplicity $m-1$. Thus, $\{q_*+1\}$
is a flat band of the operator $H$ with multiplicity $m-1$.

\emph{v}) All edges of the fundamental graph $\G_*$ with non-zero indices are loops at the
vertex $v_\n$ and their number is $2d$. The degree of this vertex $\vk_\n=\x$.
Then the number $\gz$, defined in \eqref{eq.7}, is equal to ${2d\over \x}$ and the
identity (\ref{loo}) takes the form \eqref{sp2}. \qq $\BBox$
\

\setlength{\unitlength}{0.8mm}
\begin{figure}[h]
\centering

\unitlength 1mm 
\linethickness{0.4pt}
\ifx\plotpoint\undefined\newsavebox{\plotpoint}\fi 
\begin{picture}(100,42)(0,0)
\put(-5,10){({a})}
\put(10,10){\line(1,0){30.00}} \put(10,25){\line(1,0){30.00}}
\put(10,40){\line(1,0){30.00}} \put(10,10){\line(0,1){30.00}}
\put(25,10){\line(0,1){30.00}} \put(40,10){\line(0,1){30.00}}

\put(10,10){\circle*{1}} \put(15,10){\circle*{1}}
\put(20,10){\circle*{1}} \put(25,10){\circle{1}}
\put(30,10){\circle{1}} \put(35,10){\circle{1}}
\put(40,10){\circle{1}}

\put(10,25){\circle{1}} \put(15,25){\circle{1}}
\put(20,25){\circle{1}} \put(25,25){\circle{1}}
\put(30,25){\circle{1}} \put(35,25){\circle{1}}
\put(40,25){\circle{1}}

\put(10,40){\circle{1}} \put(15,40){\circle{1}}
\put(20,40){\circle{1}} \put(25,40){\circle{1}}
\put(30,40){\circle{1}} \put(35,40){\circle{1}}
\put(40,40){\circle{1}}

\put(10,15){\circle*{1}} \put(10,20){\circle*{1}}
\put(10,30){\circle{1}} \put(10,35){\circle{1}}

\put(25,15){\circle{1}} \put(25,20){\circle{1}}
\put(25,30){\circle{1}} \put(25,35){\circle{1}}

\put(40,15){\circle{1}} \put(40,20){\circle{1}}
\put(40,30){\circle{1}} \put(40,35){\circle{1}} \put(8.0,7.0){$v_5$}
\put(13.0,7.0){$v_1$} \put(18.0,7.0){$v_2$} \put(6.0,14.0){$v_3$}
\put(6.0,19.0){$v_4$}

\put(55,10){({b})}
\put(70,10){\line(1,0){40.00}}
\put(70,9.8){\line(1,0){40.00}}
\put(70,10.2){\line(1,0){40.00}}

\put(70,9){\line(0,1){2.00}}
\put(110,9){\line(0,1){2.00}}

\put(80,10){\circle*{1.5}}
\put(100,10){\circle*{1.5}}

\put(104,12){$\s_5^0$}
\put(88,12){$\s_3^0$}

\put(72,12){$\s_1^0$}
\put(79,12){$\s_2^0$}
\put(98,12){$\s_4^0$}
\put(67,6){$\scriptstyle0$}
\put(97.3,6){$\scriptstyle3/2$}
\put(77.3,6){$\scriptstyle1/2$}
\put(109.3,6){$\scriptstyle2$}
\end{picture}
\vspace{-0.8cm} \caption{\footnotesize \rm{a}) The graph $\G$ obtained by adding
two vertices on each edge of the lattice  $\dL^2$; \quad \rm{b})
the spectrum of the Laplacian ($d=2$, $N=2$).} \label{f.10}
\end{figure}

We describe the spectrum of the Laplacian on the $d$-dimensional lattice with additional
vertices (Fig. \ref{f.10}\emph{a}).

\begin{proposition}\lb{TG2}
Let $\G$ be the graph obtained from the lattice $\dL^d$ by adding $N$ vertices on
each edge of $\dL^d$ (for $N=2$ see Fig. \ref{f.10}a). Then the
fundamental graph $\G_*$ has $\n=dN+1$ vertices and the spectrum of the Laplacian
on $\G$ has the form
\begin{equation}
\lb{acs3} \s(\D)=\s_{ac}(\D)\cup\s_{fb}(\D),
\end{equation}
where $\s_{ac}(\D)=[0,2]$ is the absolutely continuous spectrum, and the
set of all flat bands is given by{\rm :}
\begin{equation}
\lb{acs4} \s_{fb}(\D)=\Big\{1+\cos{\pi n\over  N+1}\;:
\;n=1,\ldots,N\Big\}.
\end{equation}
Here each flat band has multiplicity $d-1$ and is embedded into the
absolutely continuous spectrum.
\end{proposition}
\no \textbf{Proof.} Let $N=1$. The fundamental graph $\G_*$ has
$\n=d+1$ vertices. By (\ref{Devt}), for each
$$
\vt=(\vt_j)_{j=1}^d\in \T^d
$$
the operator $\D(\vt)$ in the standard basis \eqref{sonb} has the form
\begin{equation}\label{z2}
\D(\vt)=
\begin{pmatrix}
   1 & 0 & \ldots & -{1+e^{-i\vt_1}\over 2\sqrt{d}}\\[5pt]
   0 & 1 & \ldots & -{1+e^{-i\vt_2}\over 2\sqrt{d}}\\
   \ldots & \ldots & \ldots & \ldots\\
   -{1+e^{i\vt_1}\over 2\sqrt{d}} & -{1+e^{i\vt_2}\over 2\sqrt{d}} & \ldots & 1\\
\end{pmatrix}.
\end{equation}
Then the identity \eqref{det} yields
\begin{equation}
\begin{aligned}
\lb{sl1'} \textstyle\det\big(\D(\vt)-\l
\1_\n\big)&=(1-\l)^d\Big(1-\l-\frac1{2d(1-\l)}\,(
d+c_0)\Big)
\\
&=(1-\l)^{d-1}\Big(\l^2-2\l+{d-c_0\over 2d}\Big),
\end{aligned}
\end{equation}
where $c_0$ is defined in \eqref{fc0}. The eigenvalues of the matrix $\D(\vt)$ have
the form
\begin{equation}
\lb{ev11}
\l_1^0(\vt)=1-\sqrt{\dfrac{d+c_0}{2d}},\quad \l_{d+1}^0(\vt)=1+\sqrt{\dfrac{d+c_0}{2d}};
\end{equation}
\begin{equation}
\lb{ev2} \l_2^0(\vt)=\ldots=\l_{d}^0(\vt)=1.
\end{equation}
The identities \eqref{ev2} imply that $\{1\}$ is a flat band of the Laplacian $\D$
of multiplicity $d-1$. From the identities \eqref{ev11} we obtain that
$$
\begin{aligned}
&\l_1^{0-}=\min_{\vt\in\T^d}\l_1^0(\vt)=\l_1^0(0)=0,
\\
&\l_1^{0+}=\max_{\vt\in\T^d}\l_1^0(\vt)=\l_1^0(\vt_\pi)=1,\quad
\vt_\pi=(\pi,\ldots,\pi),
\\
&\l_{d+1}^{0-}=\min_{\vt\in\T^d}\l_{d+1}^0(\vt)=\l_{d+1}^0(\vt_\pi)=1,
\\
&\l_{d+1}^{0+}=\max_{\vt\in\T^d}\l_{d+1}^0(\vt)=\l_{d+1}^0(0)=2.
\end{aligned}
$$
Thus,
$$
\s_1^0=[0,1],\quad \s_{d+1}^0=[1,2],\quad \s_{ac}(\D)=\s_1^0\cup\s_{d+1}^0=[0,2].
$$

Let $N\geq2$. The fundamental graph $\G_*$ has $\n=dN+1$ vertices.
In this case it is more convenient to consider the operator
$J=\1-\D$ instead of the Laplacian $\D$, where $\1$ is the identity operator.

The operator $J(\vt)=\1-\D(\vt)$ in the standard basis \eqref{sonb} has the form
$$
J(\vt)=\left(
\begin{array}{cc}
  J_{dN} & y(\vt) \\
  y^*(\vt) & 0
\end{array}\right),\quad \forall\,\vt=(\vt_j)_{j=1}^d\in \T^d,
$$
where the $(dN\ts dN)$-matrix $J_{dN}=\diag(J_N,\ldots,J_N)$ and the vector $y(\vt)$
are given by
\begin{equation}
\lb{sl0}
\begin{aligned}
&J_N=\textstyle{1\over 2}\left(
\begin{array}{cccc}
  0 & 1 & 0 &\ldots \\
  1 & 0 & 1 &\ldots\\
  0 & 1 & 0 &\ldots \\
   \ldots & \ldots & \ldots &\ldots \\
\end{array}\right),\qq
y(\vt)=\ma y_1(\vt_1)\\
\vdots \\
 y_d(\vt_d) \am,
\\
&y_j(\vt_j)=\frac1{2\sqrt{d}}\ma 1\\
0\\
\vdots \\
0\\
e^{-i\vt_j} \am \in \C^N, \qq j\in\N_d.
\end{aligned}
\end{equation}
It is known that the eigenvalues of the Jacobi $({N\ts N})$-matrix $J_N$ are distinct
and have the form
\begin{equation}\label{ww.17}
\mu_n=\cos\frac{\pi n}{N+1}\in(-1,1),\qquad n=1,\ldots, N.
\end{equation}
Consequently, the matrix $J_{dN}$ has $N$ distinct eigenvalues $\m_n$ of multiplicity
$d$. Then, by the matrix property iii), the operator $J$ has at least $N=\frac{\n-1}d$
flat bands $\{\m_n\}$, $n\in\N_N$,  each of which has multiplicity $d-1$.

We describe $\s_{ac}(J)$. The identity \eqref{det} yields
\begin{equation}
\lb{sl1} \det\big(\l\1_\n-J(\vt)\big)=\big(\l-y^\ast(\vt)
(\l\1_{dN}-J_{dN})^{-1}y(\vt)\big)\det\big(\l\1_{dN}-J_{dN}\big).
\end{equation}
From the explicit form of the matrix $J_{dN}$ it follows that
\begin{equation}
\lb{sl2} \det\big(\l\1_{dN}-J_{dN}\big)=\cD^d_N(\l),\qq\textrm{
where } \ \cD_N(\l)=\det\big(\l\1_N-J_N\big).
\end{equation}
A direct calculation gives
\[
\lb{sl3} (\l\1_{dN}-J_{dN})^{-1}=\diag\,(B_N,\ldots,B_N), \quad
B_N=(\l\1_N-J_N)^{-1},
\]
\[
\lb{sl5} y^\ast(\vt) (\l\1_{dN}-J_{dN})^{-1}y(\vt)=y_1^\ast(\vt_1)
B_Ny_1(\vt_1)+\ldots+y_d^\ast(\vt_d) B_Ny_d(\vt_d),
\]
\[\lb{sl16}
y_j^\ast(\vt_j)
B_Ny_j(\vt_j)=\frac1{2^N\cD_N(\l)\,d}\big(2^{N-1}\cD_{N-1}(\l)+
\cos\vt_j\,\big),\quad j\in\N_d.
\]
Substituting \eqref{sl2}, \eqref{sl5} and \eqref{sl16} into the formula \eqref{sl1},
we obtain
\[
\lb{sl4}
\begin{aligned}
\det\big(\l\1_\n-J(\vt)\big)={\cD^{d-1}_N(\l)\over 2^N}\Big(\l\,\wt\cD_N(\l)-
\wt\cD_{N-1}(\l)- \frac{c_0}{d} \Big),
\\
 \qq\textrm{where} \qq
\wt\cD_n(\l)=2^n\cD_n(\l).
\end{aligned}
\]
The expression $\wt\cD_n(\l)$ satisfies the following recurrence
relations (the Jacobi equation):
\[
\lb{rel1}
\begin{array}{c}
\wt\cD_{n+1}(\l)=2\l\,\wt\cD_n(\l)-\wt\cD_{n-1}(\l),\quad \forall \,
 n\in\N,  \\[6pt]
\textrm{ where } \ \wt\cD_0(\l)=1,\qq \wt\cD_1(\l)=2\l.
\end{array}
\]
Thus, $\wt\cD_n(\l)$, $n=1,2,\ldots,$ are the Chebyshev polynomials
of the second kind and the following identity holds
\[
\lb{sl6} \wt\cD_n(\l)=\frac{\sin(n+1)\vp}{\sin\vp},\qq
\textrm{ where } \; \lambda=\cos\vp, \quad n\in\N.
\]
Using the formulas \eqref{rel1} and \eqref{sl6}, we rewrite \eqref{sl4} in the form
\[
\lb{sl15}
\begin{aligned}
\det\big(\l\1_\n-J(\vt)\big)&={\cD^{d-1}_N(\l)\over 2^N}
\Big(\,{1\over 2}\,\big(\wt\cD_{N+1}(\l)-
\wt\cD_{N-1}(\l)\big)- \frac{c_0}{d} \Big)
\\
& ={\cD^{d-1}_N(\l)\over 2^N}
\left(\cos(N+1)\vp-\frac{c_0}d\right).
\end{aligned}
\]
Then the eigenvalues of the matrix $J(\vt)$ are solutions of the equations
\[
\lb{spe} \cD^{d-1}_N(\l)=0\qq\textrm{ or } \qq \cos(N+1)\vp=\frac{c_0}d.
\]
From the first equation it follows that all flat bands $\m_n$, $n=1,\ldots,N$,
of the operator~$J$  are defined by the
formula~\eqref{ww.17}. Then the set of all flat bands
of the Laplacian $\D$ has the form \eqref{acs4}. Since the range of the function
$\frac{c_0}d$ is the segment $[-1,1]$, each $\l\in[-1,1]$ is a solution of the second
equation in \eqref{spe} for some $\vt\in\T^d$, i.e., an eigenvalue of the
matrix~$J(\vt)$. Thus, $\s_{ac}(J)=[-1,1]$. Consequently,
\begin{equation*}
\s_{ac}(\D)=[0,2].\qqq \BBox
\end{equation*}

\no \textbf{Proof of Theorem \ref{Tfb}.}
Items \emph{i}) -- \emph{iii}) are direct consequences of Propositions
\ref{TG1} and \ref{TG2}. \qq $\BBox$

\section{\lb{Sec7} Crystal models}
\setcounter{equation}{0}

It is known that the majority of common metals have either a face
centered cubic (FCC) structure (Fig. \ref{ff.FCC}), a body centered
cubic (BCC) structure (Fig. \ref{ff.11}), or a hexagonal close packed
(HCP) structure (see \cite{BM80}). The differences between these
structures lead to different physical properties of bulk metals. For
example, FCC metals, Cu, Au, Ag, are usually soft and ductile, which
means they can be bent and shaped easily. BCC metals are less
ductile but stronger, for example iron, while HCP metals are usually
brittle. Zinc is HCP and is difficult to bend without breaking,
unlike copper.

We describe the spectrum of the Laplace and Schr\"odinger operators on the face
centered and body centered cubic lattices.

\subsection{Body-centered cubic lattice} The body-centered cubic lattice
$\bB_2$ is obtained from the cubic lattice~$\dL^3$ by adding one vertex in the center
of each cube. This vertex is connected by an edge with each corner vertex of the cube
(Fig. \ref{ff.11}{a}). The fundamental graph of the lattice $\bB_2$ consists of
two vertices $v_1,v_2$ and 11 edges (Fig. \ref{ff.11}{b}).

\setlength{\unitlength}{1.0mm}
\begin{figure}[h]
\centering

\unitlength 1.0mm 
\linethickness{0.4pt}
\ifx\plotpoint\undefined\newsavebox{\plotpoint}\fi 
\begin{picture}(140,60)(0,0)

\put(10,10){\line(1,0){40.00}} \put(10,30){\line(1,0){40.00}}
\put(10,50){\line(1,0){40.00}}

\bezier{60}(10,10)(23.5,21.75)(37,33.5)
\bezier{60}(10,30)(23.5,21.75)(37,13.5)
\bezier{40}(17,13.5)(23.5,21.75)(30,30.0)
\bezier{40}(17,33.5)(23.5,21.75)(30,10.0)
\put(23.5,21.75){\circle{0.8}}

\bezier{60}(10,30)(23.5,41.75)(37,53.5)
\bezier{60}(10,50)(23.5,41.75)(37,33.5)
\bezier{40}(17,33.5)(23.5,41.75)(30,50.0)
\bezier{40}(17,53.5)(23.5,41.75)(30,30.0)
\put(23.5,41.75){\circle{0.8}}

\bezier{60}(30,10)(43.5,21.75)(57,33.5)
\bezier{60}(30,30)(43.5,21.75)(57,13.5)
\bezier{40}(37,13.5)(43.5,21.75)(50,30.0)
\bezier{40}(37,33.5)(43.5,21.75)(50,10.0)
\put(43.5,21.75){\circle*{1}}

\bezier{60}(30,30)(43.5,41.75)(57,53.5)
\bezier{60}(30,50)(43.5,41.75)(57,33.5)
\bezier{40}(37,33.5)(43.5,41.75)(50,50.0)
\bezier{40}(37,53.5)(43.5,41.75)(50,30.0)
\put(43.5,41.75){\circle{0.8}}

\bezier{60}(17,33.5)(37,33.5)(57,33.5)
\bezier{60}(17,13.5)(37,13.5)(57,13.5)

\put(17,53.5){\line(1,0){40.00}} \put(10,10){\line(0,1){40.00}}
\put(30,10){\line(0,1){40.00}} \put(50,10){\line(0,1){40.00}}
\put(57,13.5){\line(0,1){40.00}}
\bezier{12}(10,10)(13.5,11.75)(17,13.5)

\bezier{60}(17,13.5)(17,33.5)(17,53.5)
\bezier{60}(37,13.5)(37,33.5)(37,53.5)

\bezier{12}(30,10)(33.5,11.75)(37,13.5)
\put(50,10){\line(2,1){7.00}}
\bezier{12}(10,30)(13.5,31.75)(17,33.5)
\bezier{12}(30,30)(33.5,31.75)(37,33.5)

\put(50,30){\line(2,1){7.00}} \put(10,50){\line(2,1){7.00}}
\put(30,50){\line(2,1){7.00}} \put(50,50){\line(2,1){7.00}}
\put(10,10){\circle{0.8}} \put(30,10){\circle{0.8}}
\put(50,10){\circle*{1}} \put(10,30){\circle{0.8}}
\put(30,30){\circle{0.8}} \put(50,30){\circle{0.8}}
\put(10,50){\circle{0.8}} \put(30,50){\circle{0.8}}
\put(50,50){\circle{0.8}}

\put(17,53.5){\circle{0.8}} \put(37,53.5){\circle{0.8}}
\put(57,53.5){\circle{0.8}} \put(17,33.5){\circle{0.8}}
\put(37,33.5){\circle{0.8}} \put(57,33.5){\circle{0.8}}
\put(17,13.5){\circle{0.8}} \put(37,13.5){\circle{0.8}}
\put(57,13.5){\circle{0.8}}

\put(50,10){\vector(0,1){20.00}} \put(50,10){\vector(-1,0){20.00}}
\put(50,10){\vector(2,1){7.00}}

\put(54,10.0){$\scriptstyle a_1$}
\put(37.0,8.0){$\scriptstyle a_2$}
\put(51.0,24){$\scriptstyle a_3$}
\put(42.0,18.5){$\scriptstyle v_1$}

\put(49.0,7.5){$\scriptstyle v_2$}

\put(2,5){({a})}
\bezier{160}(80,10)(93.5,21.75)(107,33.5)
\bezier{160}(80,30)(93.5,21.75)(107,13.5)
\bezier{140}(87,13.5)(93.5,21.75)(100,30.0)
\bezier{140}(87,33.5)(93.5,21.75)(100,10.0)
\put(93.5,21.75){\circle*{1.0}}

\bezier{30}(87,13.5)(97,13.5)(107,13.5)
\bezier{30}(87,13.5)(87,23.5)(87,33.5)
\bezier{12}(87,13.5)(83.5,11.75)(80,10)

\put(100,10){\vector(2,1){7.00}}
\bezier{12}(100,30)(103.5,31.75)(107,33.5)
\bezier{12}(80,30)(83.5,31.75)(87,33.5)
\bezier{30}(80,10)(80,20)(80,30)
\bezier{30}(107,13.5)(107,23.5)(107,33.5)
\put(100,10){\vector(-1,0){20.00}}
\bezier{30}(87,33.5)(97,33.5)(107,33.5)

\put(100,10){\vector(0,1){20.00}} \bezier{30}(80,30)(90,30)(100,30)
\put(80,10){\circle{1}} \put(80,30){\circle{1}}
\put(100,30){\circle{1}} \put(87,13.5){\circle{1}}
\put(107,13.5){\circle{1}} \put(107,33.5){\circle{1}}
\put(87,33.5){\circle{1}}

\put(100,10){\circle*{1}}

\put(100.0,7.5){$\scriptstyle v_2$} \put(92.1,18.0){$\scriptstyle
v_1$}

\put(103.0,9.5){$\scriptstyle a_1$} \put(89.0,8.0){$\scriptstyle
a_2$} \put(101.0,21.0){$\scriptstyle a_3$}

\put(70,5){({b})}
\end{picture}

\caption{\footnotesize {a}) The body-centered cubic lattice $\bB_2$; {b})
  the fundamental graph.} \label{ff.11}
\end{figure}

\begin{proposition}\lb{BCC}
i) The spectrum of the Laplacian on the body-centered cubic lattice $\bB_2$
has the form
\[
\lb{sBCC}\textstyle
\s(\D)=\s_{ac}(\D)=\s_1^0\cup\s_2^0=[0,{11\over 7}],\quad  \s_1^0=[0,1],
\quad \s_2^0=[1,{11\over 7}].
\]
ii) Let the Schr\"odinger operator $H=\D+Q$ act on the
body-centered cubic lattice $\bB_2$, and let the potential $Q$
satisfy
$$
Q(v_j)=q_j, \quad j\in \N_2, \quad q_2=0.
$$
Then the spectrum of the operator $H$  has the form
\[
\lb{SB} \s(H)=\s_{ac}(H)=\s_1\cup\s_2=[\l^-_1,\l^+_2]\sm \g_1,
\]
where
$$
\begin{aligned}
 \l_1^-&=\frac{11}{14}+\frac12\,q_1-\frac12
\sqrt{\Big(\frac37+q_1\Big)^2+ \frac{16}7},
\\
\l_2^+&=\frac{11}{14}+\frac12\,q_1+\frac12
\sqrt{\Big(\frac37+q_1\Big)^2+ \frac{16}7},
\end{aligned}
$$
and the gap $\g_1$ is given by
$$
\g_1=(\l^+_1,\l^-_2)=\ca ({10\over 7},1+q_1), & \textrm{ if } \qq {3\over 7}<q_1\\[2pt]
(1+q_1,1+q_1)=\varnothing , & \textrm{ if } \qq -{1\over 7}\leq q_1\leq{3\over 7}\\[2pt]
(1+q_1,{6\over 7}),& \textrm{ if } \qq q_1<-{1\over 7}  \ac.
$$
\end{proposition}

\no \textbf{Proof.} \emph{i}) The fundamental graph $\G_*$ of the body-centered cubic lattice
$\bB_2$ consists of two vertices $v_1,v_2$ with degrees $\vk_1=8$, $\vk_2=14$; 11
oriented edges
$$
\be_1=\be_2=\be_3=(v_2,v_2),\quad \be_4=\ldots=\be_{11}=(v_1,v_2)
$$
and their inverse edges (Fig. \ref{ff.11}\emph{b}). The indices of
the fundamental graph edges are given by
$$
\begin{array}{ll}
\t(\be_1)=\t(\be_5)=(1,0,0),\qquad & \t(\be_2)=\t(\be_6)=(0,1,0),\\[4pt]
\t(\be_3)= \t(\be_7)=(0,0,1), \quad & \t(\be_4)=(0,0,0),\\[4pt]
\t(\be_8)=(1,1,0),\qquad\qquad\;\; & \t(\be_9)=(1,0,1),\\[4pt]
\t(\be_{10})=(0,1,1),\qquad\qquad\, & \t(\be_{11})=(1,1,1).
\end{array}
$$
For each  $\vt=(\vt_j)_{j=1}^3\in\T^3$ the operator $\D(\vt)$ in the standard basis
\eqref{sonb} has the form
\[
\lb{el}
\begin{aligned}
& \D(\vt)=\left(
\begin{array}{cc}
1  & \D_{12} \\[6pt]
\bar\D_{12} & \D_{22}
\end{array}\right)(\vt),
\\
&    \ca
\D_{12}(\vt)=-\frac1{4\sqrt{7}}\,(1+e^{-i\vt_1})(1+e^{-i\vt_2})(1+e^{-i\vt_3})\\
\D_{22}(\vt)=1-\frac{c_0}{7}\,\ac,
\end{aligned}
\]
where
\[
\lb{pol} c_0=c_1+c_2+c_3,\quad c_j=\cos\vt_j,\quad j=1,2,3.
\]
By a direct calculation, we get
$$\textstyle
\det(\D(\vt)-\l
\1_2)=(\l-1)^2+{c_0\over 7}(\l-1)-\frac1{14}\,(1+c_1)(1+c_2)(1+c_3).
$$
The eigenvalues of each matrix $\D(\vt)$ are given by
\[
\lb{ev}\textstyle
\l_n^0(\vt)=1-{c_0\over 14}+{(-1)^n\over 2}\sqrt{{c_0^2\over 49}+
{2\over 7}\,(1+c_1)(1+c_2)(1+c_3)}\;,  \quad n=1,2.
\]
Thus, the spectrum of the Laplacian on the body-centered cubic
lattice $\bB_2$ has the form
$$
\begin{aligned}
&\s(\D)=\s_{ac}(\D)=[\l_1^{0-},\l_1^{0+}]\cup[\l_2^{0-},\l_2^{0+}],
\\
&[\l_n^{0-},\l_n^{0+}]=\l_n^0(\T^3),\quad n=1,2.
\end{aligned}
$$
From the matrix property iii) it follows that
\[
\lb{razz} \l_1^0(\vt)\le1\le\l_2^0(\vt),\quad \forall\,\vt\in\T^3.
\]
Investigating the function $\l_2^0$ for the maximum and using
\eqref{bes}, \eqref{razz}, we obtain
$$
\begin{array}{ll}
\l_1^{0-}=\l_1^0(0)=0, & \l_1^{0+}=\l_1^0(\pi,\pi,\pi)=1,\\[4pt]
\l_2^{0-}=\l_2^0(\pi,0,0)=1,\qquad
&\l_2^{0+}=\l_2^0(0,0,0)={11\over 7}.
\end{array}
$$

\emph{ii}) For each $\vt\in\T^3$ the operator $H(\vt)$ in the standard basis \eqref{sonb}
has the form
$$
H(\vt)=
\begin{pmatrix}
1+q_1  & \D_{12} \\[6pt]
\bar\D_{12} & \D_{22}
\end{pmatrix}(\vt),
$$
where $\D_{12}(\vt),\D_{22}(\vt)$ are given in \eqref{el}. By a direct
calculation, we get
$$
\det(H(\vt)-\l \1_2)=(\l-1)^2+\Big({c_0\over 7}-q_1\Big)(\l-1)
-{c_0\over 7}q_1-\frac1{14}\,(1+c_1)(1+c_2)(1+c_3),
$$
where $c_0$ and $c_j$, $j=1,2,3$, are defined by \eqref{pol}.
The eigenvalues of each matrix $H(\vt)$ are given by
\begin{equation}
\lb{ev1}
\l_n(\vt)=1-\frac{c_0}{14}+\frac12\,q_1+\frac{(-1)^n}2
\sqrt{\Big({c_0\over 7}+q_1\Big)^2+ \frac27\,(1+c_1)(1+c_2)(1+c_3)},
\qq n=1,2.
\end{equation}
By \eqref{bes} and the matrix property iii), we have
$$
\begin{aligned}
&\l_1^-=\l_1(0)=\frac{11}{14}+\frac12\,q_1-\frac12
\sqrt{\Big(\frac37+q_1\Big)^2+ \frac{16}7}\;,
\\
&\l_1(\vt)\leq1+q_1\leq\l_2(\vt),\quad \forall\,\vt\in\T^3.
\end{aligned}
$$
Then an investigation of the functions $\l_1$, $\l_2$ for the extremes
gives
$$
\l_1^+=\l_1(\pi,\pi,\pi)=
\left\{\begin{array}{rl}
         1+q_1, & \textrm{if} \qq q_1\leq{3\over 7}\,\\[4pt]
         {10\over 7}, & \textrm{if} \qq  q_1>{3\over 7}\,
       \end{array}\right.,
$$
$$
\l_2^-=\l_2(\pi,0,0)=
\left\{\begin{array}{rl}
         1+q_1, & \textrm{if} \qq  q_1\geq -{1\over 7}\,\\[4pt]
         {6\over 7}, & \textrm{if} \qq  q_1<-{1\over 7}
       \end{array}\right.,
$$
$$
\l_2^+=\textstyle\l_2(0)=
\frac{11}{14}+\frac12\,q_1+\frac12 \sqrt{\big(\frac37+q_1\big)^2+
\frac{16}7}\;.
$$
This proves item ii). \qq $\BBox$

\subsection{Face-centered cubic lattice} The face-centered cubic lattice $\bB_4$
is obtained from the cubic lattice $\dL^3$ by adding one vertex at the center of each cube
face. This vertex is connected by an edge with each corner vertex of the cube
face (Fig.~\ref{ff.FCC}{a}). The fundamental graph  of the lattice $\bB_4$ consists
of four vertices $v_1,v_2,v_3,v_4$ and 15 edges (Fig. \ref{ff.FCC}{b}).

\setlength{\unitlength}{1.0mm}
\begin{figure}[h]
\centering

\unitlength 1.0mm 
\linethickness{0.4pt}
\ifx\plotpoint\undefined\newsavebox{\plotpoint}\fi 
\begin{picture}(140,60)(0,0)

\bezier{80}(17,13.5)(37,33.5)(57,53.5)
\bezier{40}(17,33.5)(27,43.5)(37,53.5)
\bezier{40}(37,13.5)(47,23.5)(57,33.5)

\bezier{80}(57,13.5)(37,33.5)(17,53.5)
\bezier{40}(17,33.5)(27,23.5)(37,13.5)
\bezier{40}(37,53.5)(47,43.5)(57,33.5)

\bezier{40}(10,10)(23.5,11.75)(37,13.5)
\bezier{40}(30,10)(43.5,11.75)(57,13.5)
\bezier{20}(30,10)(23.5,11.75)(17,13.5)
\bezier{20}(50,10)(43.5,11.75)(37,13.5)

\bezier{40}(10,30)(23.5,31.75)(37,33.5)
\bezier{40}(30,30)(43.5,31.75)(57,33.5)
\bezier{20}(30,30)(23.5,31.75)(17,33.5)
\bezier{20}(50,30)(43.5,31.75)(37,33.5)

\bezier{180}(10,50)(23.5,51.75)(37,53.5)
\bezier{180}(30,50)(43.5,51.75)(57,53.5)
\bezier{180}(30,50)(23.5,51.75)(17,53.5)
\bezier{180}(50,50)(43.5,51.75)(37,53.5)

\bezier{40}(10,10)(13.5,21.75)(17,33.5)
\bezier{40}(10,30)(13.5,21.75)(17,13.5)
\bezier{40}(10,30)(13.5,41.75)(17,53.5)
\bezier{40}(10,50)(13.5,41.75)(17,33.5)

\bezier{40}(30,10)(33.5,21.75)(37,33.5)
\bezier{40}(30,30)(33.5,21.75)(37,13.5)
\bezier{40}(30,30)(33.5,41.75)(37,53.5)
\bezier{40}(30,50)(33.5,41.75)(37,33.5)

\bezier{140}(50,10)(53.5,21.75)(57,33.5)
\bezier{140}(50,30)(53.5,21.75)(57,13.5)
\bezier{140}(50,30)(53.5,41.75)(57,53.5)
\bezier{140}(50,50)(53.5,41.75)(57,33.5)

\put(10,10){\line(1,1){40.00}} \put(10,30){\line(1,-1){20.00}}
\put(30,50){\line(1,-1){20.00}} \put(10,50){\line(1,-1){40.00}}
\put(10,30){\line(1,1){20.00}} \put(30,10){\line(1,1){20.00}}

\put(10,10){\line(1,0){40.00}} \put(10,30){\line(1,0){40.00}}
\put(10,50){\line(1,0){40.00}}
\bezier{60}(17,33.5)(37,33.5)(57,33.5)
\bezier{60}(17,13.5)(37,13.5)(57,13.5)
\put(17,53.5){\line(1,0){40.00}}
\put(10,10){\line(0,1){40.00}} \put(30,10){\line(0,1){40.00}}
\put(50,10){\line(0,1){40.00}} \put(57,13.5){\line(0,1){40.00}}
\bezier{12}(10,10)(13.5,11.75)(17,13.5)
\bezier{60}(17,13.5)(17,33.5)(17,53.5)
\bezier{60}(37,13.5)(37,33.5)(37,53.5)

\bezier{12}(30,10)(33.5,11.75)(37,13.5)
\put(50,10){\line(2,1){7.00}}
\bezier{12}(10,30)(13.5,31.75)(17,33.5)
\bezier{12}(30,30)(33.5,31.75)(37,33.5)

\put(50,30){\line(2,1){7.00}}
\put(10,50){\line(2,1){7.00}} \put(30,50){\line(2,1){7.00}}
\put(50,50){\line(2,1){7.00}}
\put(10,10){\circle{0.8}} \put(30,10){\circle{0.8}}
\put(50,10){\circle{0.8}}
\put(10,30){\circle{0.8}} \put(30,30){\circle{0.8}}
\put(50,30){\circle{0.8}}
\put(10,50){\circle{0.8}} \put(30,50){\circle{0.8}}
\put(50,50){\circle{0.8}}
\put(20,20){\circle{0.8}} \put(40,20){\circle{0.8}}

\put(20,40){\circle{0.8}} \put(40,40){\circle{0.8}}

\put(27,23.5){\circle{0.8}} \put(47,23.5){\circle{0.8}}

\put(27,43.5){\circle{0.8}} \put(47,43.5){\circle{0.8}}

\put(23.5,11.75){\circle{0.8}} \put(43.5,11.75){\circle{0.8}}
\put(23.5,31.75){\circle{0.8}} \put(43.5,31.75){\circle{0.8}}
\put(23.5,51.75){\circle{0.8}} \put(43.5,51.75){\circle{0.8}}

\put(13.5,21.75){\circle{0.8}} \put(13.5,41.75){\circle{0.8}}
\put(33.5,21.75){\circle{0.8}} \put(33.5,41.75){\circle{0.8}}
\put(53.5,21.75){\circle{0.8}} \put(53.5,41.75){\circle{0.8}}

\put(17,53.5){\circle{0.8}} \put(37,53.5){\circle{0.8}}
\put(57,53.5){\circle{0.8}} \put(17,33.5){\circle{0.8}}
\put(37,33.5){\circle{0.8}} \put(57,33.5){\circle{0.8}}
\put(17,13.5){\circle{0.8}} \put(37,13.5){\circle{0.8}}
\put(57,13.5){\circle{0.8}}

\put(50,10){\vector(0,1){20.00}} \put(50,10){\vector(-1,0){20.00}}
\put(50,10){\vector(2,1){7.00}}

\put(49,7.0){$\scriptstyle v_4$} \put(53,10.0){$\scriptstyle a_1$}
\put(35.0,7.5){$\scriptstyle a_2$} \put(46.5,19){$\scriptstyle
a_3$}

\put(50,10){\circle*{1}} \put(40,20){\circle*{1}}
\put(42.5,13.0){$\scriptstyle v_1$}
\put(38.7,21.7){$\scriptstyle v_3$}
\put(54.2,21.3){$\scriptstyle v_2$}
\put(43.5,11.75){\circle*{1}}
\put(53.5,21.75){\circle*{1}}

\put(2,5){({a})}
\put(80,10){\line(1,1){20.00}} \put(100,10){\line(-1,1){20.00}}
\bezier{180}(80,10)(93.5,11.75)(107,13.5)
\bezier{180}(87,13.5)(93.5,11.75)(100,10)
\bezier{30}(87,13.5)(97,13.5)(107,13.5)
\bezier{30}(87,13.5)(87,23.5)(87,33.5)
\bezier{12}(87,13.5)(83.5,11.75)(80,10)

\put(100,10){\vector(2,1){7.00}}
\bezier{12}(100,30)(103.5,31.75)(107,33.5)
\bezier{12}(80,30)(83.5,31.75)(87,33.5)
\bezier{30}(80,10)(80,20)(80,30)
\bezier{30}(107,13.5)(107,23.5)(107,33.5)
\put(100,10){\vector(-1,0){20.00}}
\bezier{30}(87,33.5)(97,33.5)(107,33.5)

\bezier{180}(100,10)(103.5,21.75)(107,33.5)
\bezier{180}(100,30)(103.5,21.75)(107,13.5)
\put(100,10){\vector(0,1){20.00}} \bezier{30}(80,30)(90,30)(100,30)
\put(80,10){\circle{1}} \put(80,30){\circle{1}}
\put(100,30){\circle{1}} \put(87,13.5){\circle{1}}
\put(107,13.5){\circle{1}} \put(107,33.5){\circle{1}}
\put(87,33.5){\circle{1}}
\put(100,10){\circle*{1}} \put(90,20){\circle*{1}}
\put(100.0,8.0){$\scriptstyle v_4$} \put(92.5,13.0){$\scriptstyle
v_1$} \put(88.7,21.7){$\scriptstyle v_3$}
\put(104.2,21.3){$\scriptstyle v_2$}
\put(104.0,10.5){$\scriptstyle a_1$} \put(85.0,8.0){$\scriptstyle
a_2$} \put(97.0,25.0){$\scriptstyle a_3$}
\put(93.4,11.8){\circle*{1}} \put(103.5,21.8){\circle*{1}}

\put(70,5){({b})}
\end{picture}

\vspace{-0.5cm} \caption{\footnotesize {a}) The face-centered cubic lattice;\quad
{b}) the fundamental graph.} \label{ff.FCC}
\end{figure}

\begin{proposition}\lb{FCC}
i) The spectrum of the Laplacian $\D$ on the face-centered cubic
lattice $\bB_4$ has the form
\begin{equation}
\begin{aligned}
\lb{sBCC} 
&\s(\D)=\s_{ac}(\D)\cup\s_{fb}(\D),
\\
&\s_{ac}(\D)=[0,1]\cup\bigg[\frac{4}{3},\frac{5}{3}\bigg],
\qqq \s_{fb}(\D)=\{1\},
\end{aligned}
\end{equation}
where the flat band $\{1\}$ has multiplicity $2$.

ii) Let the Schr\"odinger operator $H=\D+Q$ act on the
face-centered cubic lattice $\bB_4$, and let the potential $Q$
satisfy
$$
Q(v_j)=q_j, \quad j\in \N_4, \quad q_4=0.
$$
Then the spectrum of the operator $H$ has the form
\[
\lb{sBCC1}\textstyle \s(H)=\s_{ac}(H)\cup\s_{fb}(H),
\]
where
\begin{equation}
\lb{fla}
\s_{fb}(H)
=  \begin{cases}
    \varnothing , & \text{if }\  q_1,q_2,q_3 \text{ are distinct (a generic potential);}\\ 
    \{q_{(1)}+1\}, & \text{if }q_j=q_k=q_{(1)}\neq q_n \text{ for some } j,k,n=1,2,3,\\   
    & n\ne j\neq k, \qq k\neq n; \\
    \{q_{(2)}+1\} , & \textrm{if }\  q_1=q_2=q_3=q_{(2)}.  \\
  \end{cases}
\end{equation}
Moreover, the flat bands $\{q_{(1)}+1\}$ and $\{q_{(2)}+1\}$ have
multiplicities $1$ and $2$, respectively.
\end{proposition}

\no \textbf{Proof.} \emph{i}) The fundamental graph $\G_*$ of the face-centered cubic lattice
$\bB_4$ consists of four vertices $v_1,v_2,v_3,v_4$ with degrees $\vk_1=\vk_2=\vk_3=4$,
$\vk_4=18$; 15 oriented edges
\begin{equation*}
\begin{aligned}
\be_1&=\be_2=\be_3=(v_4,v_4),\quad \be_4=\ldots=\be_7=(v_1,v_4),
\\
\be_8&=\ldots=\be_{11}=(v_2,v_4),\quad
\be_{12}=\ldots=\be_{15}=(v_3,v_4)
\end{aligned}
\end{equation*}
and their inverse edges (Fig. \ref{ff.FCC}{b}). The indices of
the fundamental graph edges are given by:
$$
\begin{aligned}
&\t(\be_1)=\t(\be_5)=\t(\be_9)=(1,0,0),&\qquad \t(\be_7)=(1,1,0),
\\
& \t(\be_2)=\t(\be_6)=
\t(\be_{13})=(0,1,0),&\qquad  \t(\be_{11})=(1,0,1),\\
&\t(\be_3)=\t(\be_{10})=\t(\be_{14})=(0,0,1),&\qquad\t(\be_{15})=(0,1,1).
\\
&\t(\be_4)=\t(\be_8)=\t(\be_{12})=(0,0,0),
\end{aligned}
$$
For each $\vt=(\vt_j)_{j=1}^3\in\T^3$ the operator $\D(\vt)$ in the standard basis
\eqref{sonb} has the form
\begin{equation}
\begin{aligned}
\lb{DDD}
&\D(\vt)=
\begin{pmatrix}
1 & 0  & 0 & \D_{14} \\[6pt]
0 & 1  & 0 & \D_{24} \\[6pt]
0 & 0  & 1 & \D_{34} \\[6pt]
\bar\D_{14} & \bar\D_{24} & \bar\D_{34} & \D_{44} \\[6pt]
\end{pmatrix}(\vt),
\\
&
\ca \D_{14}(\vt)=-\frac1{6\sqrt{2}}\,(1+e^{-i\vt_1})(1+e^{-i\vt_2}),\\
\D_{24}(\vt)=-\frac1{6\sqrt{2}}\,(1+e^{-i\vt_1})(1+e^{-i\vt_3}),\\
\D_{34}(\vt)=-\frac1{6\sqrt{2}}\,(1+e^{-i\vt_2})(1+e^{-i\vt_3}),\\
\D_{44}(\vt)=1-\frac19\,c_0,\ac
\end{aligned}
\end{equation}
where $c_0=c_1+c_2+c_3$, $c_j=\cos\vt_j$,
$j=1,2,3$. By a direct calculation, we get
$$
\det(\D(\vt)-\l
\1_4)=(\l-1)^2\Big((\l-1)^2+{c_0\over 9}(\l-1)-\frac16-{c_0\over 9}-
\frac1{18}\eta\Big),
$$
where
$$
\eta=\eta(\vt)=c_1c_2+c_1c_3+c_2c_3.
$$
From the matrix property iii) it follows that
\[
\lb{raz2} \l_1^0(\vt)\le1\le\l_2^0(\vt)\le1\le\l_3^0(\vt)\le1\le\l_4^0(\vt),\quad \forall\,\vt\in\T^3.
\]
Then the eigenvalues of the matrix $\D(\vt)$ are given by
\begin{equation}
\lb{evf}
\begin{aligned}
\l_{1+3s}^0(\vt)&=1-{c_0\over 18}+{(-1)^{s+1}\over 2}\sqrt{\Big({c_0\over 9}\Big)^2+{4c_0\over 9}+{2\over 3}+
{2\over 9}\eta}\;,
\\
 s&=0,1,\quad \l_2^0(\vt)=\l_3^0(\vt)=1.
\end{aligned}
\end{equation}
Thus, the spectrum of the Laplacian on the face-centered cubic
lattice $\bB_4$ has the form
$$
\begin{aligned}
&\s(\D)=\s_{ac}(\D)\cup\s_{fb}(\D),
\\
&\s_{ac}(\D)=[\l_1^{0-},\l_1^{0+}]\cup[\l_4^{0-},\l_4^{0+}],
\qq
\s_{fb}(\D)=\{1\},
\end{aligned}
$$
where the flat band $\{1\}$ has multiplicity~2 and
$$
[\l_n^{0-},\l_n^{0+}]=\l_n^0(\T^3), \quad n=1,4.
$$
Next, using \eqref{bes} and \eqref{raz2}, we obtain
$$
\l_1^{0-}=\l_1^0(0)=0,\quad \l_1^{0+}=\l_1^0(\pi,\pi,\pi)=1.
$$
A direct calculation gives
$$\textstyle
\l_4^{0-}=\l_4^0(\pi,\pi,\pi)={4\over 3},\quad
\l_4^{0+}=\l_4^0(0)={5\over 3}\,.
$$

\emph{ii}) For each $\vt\in\T^3$ the operator $H(\vt)$ has the form
$$
H(\vt)=\D(\vt)+q,\quad q=\diag\,(q_1,q_2,q_3,0),
$$
where $\D(\vt)$ is defined by \eqref{DDD}. We write the characteristic
polynomial of the matrix $H(\vt)$ in the form of a linear
combination of linearly independent functions:
\begin{equation}\lb{lco}
\begin{aligned}
\det(H(\vt)-\l\1_4)&=\eta_1(\l,q)+\eta_2(\l,q)c_1+\eta_3(\l,q)c_2+
\eta_4(\l,q)c_3
\\
&+\eta_5(\l,q)c_1c_2+ \eta_6(\l,q)c_1c_3+\eta_7(\l,q)c_2c_3,
\end{aligned}
\end{equation}
where $c_j=\cos\vt_j$, $j=1,2,3$,
\begin{equation*}
\begin{aligned}
 \eta_1(\l,q)&=(\l-1)^4\!-\!\z_1(\l-1)^3\! +\!\Big(\z_2\!-\!{1\over 6}\Big)(\l\!-\!1)^2
\!+\!\Big({1\over 9}\z_1-\z_3\Big)(\l-1)-{1\over 18}\z_2,
\\
 \eta_2(\l,q)&={1\over 9}(\l-1)^3-{1\over 9}(\z_1+1)(\l-1)^2
 \\
 &\quad+
\Big({1\over 18}\,\z_1+{1\over 9}\,\z_2+{1\over 18}\,q_3\Big)(\l-1)
-{1\over 18}\,\z_2-{1\over 9}\,\z_3+{1\over 18}\,q_1q_2,
\\
\eta_3(\l,q)&={1\over 9}(\l-1)^3-{1\over 9}(\z_1+1)(\l-1)^2
\\
&\quad+
\Big({1\over 18}\,\z_1+{1\over 9}\,\z_2+{1\over 18}\,q_2\Big)(\l-1)
-{1\over 18}\,\z_2-{1\over 9}\,\z_3+{1\over 18}\,q_1q_3,
\\
 \eta_4(\l,q)&={1\over 9}(\l-1)^3-{1\over 9}(\z_1+1)(\l-1)^2
 \\
 &\quad+
\Big({1\over 18}\,\z_1+{1\over 9}\,\z_2+{1\over 18}\,q_1\Big)(\l-1)
-{1\over 18}\,\z_2-{1\over 9}\,\z_3+{1\over 18}\,q_2q_3,
\\
\eta_5(\l,q)&=-{1\over 18}\,(\l-1)^2+{1\over 18}\,(q_2+q_3)(\l-1)-{1\over 18}\,q_2q_3,
\\
\eta_6(\l,q)&=-{1\over 18}\,(\l-1)^2+{1\over 18}\,(q_1+q_3)(\l-1)-{1\over 18}\,q_1q_3,
\\
\eta_7(\l,q)&=-{1\over 18}\,(\l-1)^2+{1\over 18}\,(q_1+q_2)(\l-1)-{1\over 18}\,q_1q_2,
\\
\z_1&=q_1+q_2+q_3,\\
\z_2&=q_1q_2+q_1q_3+q_2q_3,\\
\z_3&=q_1q_2q_3.
\end{aligned}
\end{equation*}

A point $\l$ is a flat band of the operator $H$ if and only if
$$
\det\,(H(\vt)-\l\1_4)=0
$$
for all $\vt\in\T^3$. Since the linear combination \eqref{lco} of linearly
independent functions is equal to 0, we obtain the system of equations
$$
\eta_s(\l,q)=0,\quad s\in\N_7.
$$
 All solutions of this system have the form
$$
\l=\left\{
  \begin{array}{cl}
    q_{(1)}+1, & \textrm{ if $q_j=q_k=q_{(1)}\neq q_n$ for some $j,k,n=1,2,3$, }\\
    & \; n\ne j\neq k, \qq k\neq n; \\[4pt]
    q_{(2)}+1 , & \textrm{ if $q_1=q_2=q_3=q_{(2)}$.  }\\
  \end{array}\right.
$$
The roots $q_{(1)}+1$ and $q_{(2)}+1$ have multiplicities 1 and 2, respectively.
This proves item ii). \qq $\BBox$

\section{\lb{Sec8} Properties of  matrices}
\setcounter{equation}{0}
We denote by
$$
\l_1(A)\leq\ldots\leq\l_\n(A)
$$
the eigenvalues of a self-adjoint ($\n\ts\n$)-matrix $A$, arranged in
non-decreasing order, counting multiplicities. The following well-known properties
of matrices hold.

\bigskip

{i) {\it For each $n\in\N_\n$ the eigenvalue $\l_n(A)$ satisfies the
minimax principle}:
\begin{equation}
\lb{CF1} \l_n(A)=\min_{S_n\subset\C^\n}\max_{\|x\|=1 \atop x\in
S_n}\lan Ax,x\rangle,
\end{equation}
\begin{equation}
\lb{CF2} \l_n(A)=\max_{S_{\n-n+1}\subset\C^\n}\min_{\|x\|=1 \atop
x\in S_{\n-n+1}}\lan Ax,x\rangle,
\end{equation}
{\it where $S_n$ denotes a subspace of dimension $n$ and the outer
optimization is over all subspaces of the indicated dimension}
({\it see }\cite[{\it p.}~180]{HJ85}).

\bigskip

{ii) {\it  Let $A,B$ be self-adjoint $(\nu\ts\nu)$-matrices. Then for each
$n\in\N_\n$ we have}
$$
\l_n(A)+\l_1(B)\leq\l_n(A+B)\leq\l_n(A)+\l_\n(B)
$$
({\it see}~\cite[{\it Theorem}~4.3.1]{HJ85}).

\bigskip

{iii)  {\it Let $B=\ma   A & y \\
  y^\ast & a \am$ be a self-adjoint $(\nu+1)\ts(\nu+1)$-matrix for some
  self-adjoint $(\nu\ts \nu)$-matrix $A$, some real number $a$
and some vector $y\in\C^{\nu}$. Then}
$$
\l_1(B)\leq\l_1(A)\leq\l_2(B)\leq\ldots\leq\l_\n(B)\leq\l_\n(A)\leq\l_{\nu+1}(B)
$$
({\it see} \cite[{\it Theorem}~4.3.8]{HJ85}).

{iv) {\it Let $A=\{A_{jk}\}$ be a self-adjoint $(\nu\ts\nu)$-matrix, and let
$\G(A)$ be a graph on $\n$ vertices $v_1,\ldots,v_\n$ such that there is an edge
$(v_j,v_k)$ in it if and only if $A_{jk}\neq0$. Then $A$ is irreducible if and only if
the graph $\G(A)$ is connected} ({\it see}~\cite[{\it Theorem}~6.2.24]{HJ85}).

{v) {\it Let $A$ be an irreducible self-adjoint $(\nu\ts\nu)$-matrix with nonnegative
entries. Then the largest eigenvalue $\l=\|A\|$ of the matrix $A$ is simple,
and for some vector $x$ with positive components we have} $Ax=\l x$ ({\it see }
\cite[{\it Theorem}~8.4.4]{HJ85}).

{vi) {\it Let $M=\ma   A & B \\
         C & D       \am$ be a $(\n\ts\n)$-matrix
for some square matrices $A,D$ and some matrices $B,C$. Then}
\[
\lb{det} \det M=\det A\cdot\det\big(D-CA^{-1}B\big)
\]
({\it see}~\cite[{\it pp.}~21--22]{HJ85}).

{vii) {\it Let $V=\{V_{jk}\}$ be a self-adjoint $(\n\ts \n)$-matrix, and let
$$
B=\diag\, \{B_1,\ldots, B_\n\},\quad B_j=\sum\limits_{k=1}^\n |V_{jk}|.
$$
Then the following estimates hold}:
$$
-B\le V\le B.
$$
 ({\it see, e.g.,}~\cite{K13}).

\medskip

\footnotesize
\textbf{Acknowledgments.}  Our study was
supported by the RSF grant  No. 18-11-00032.

\end{document}